\documentclass[11pt]{article}
\usepackage{amsfonts}
\usepackage{amsfonts}
\usepackage{amsfonts}
\usepackage{mathrsfs}
\usepackage{bm}
\usepackage{amssymb,amsmath} 
\usepackage{cite}
 \allowdisplaybreaks

\catcode`!=11
\let\!int\int \def\int{\displaystyle\!int}
\let\!lim\lim \def\lim{\displaystyle\!lim}
\let\!sum\sum \def\sum{\displaystyle\!sum}
\let\!sup\sup \def\sup{\displaystyle\!sup}
\let\!inf\inf \def\inf{\displaystyle\!inf}
\let\!cap\cap \def\cap{\displaystyle\!cap}
\let\!max\max \def\max{\displaystyle\!max}
\let\!min\min \def\min{\displaystyle\!min}
\let\!frac\frac \def\frac{\displaystyle\!frac}
\catcode`!=12

\let\oldsection\section
\renewcommand\section{\setcounter{equation}{0}\oldsection}

\allowdisplaybreaks
\def\pf{\it{Proof.}\rm\quad}

\def\R{\mathbb{R}}
\def\N{\mathbb{N}}

\newtheorem{thm}{Theorem}[section]
\newtheorem{lem}[thm]{Lemma}
\newtheorem{cor}[thm]{Corollary}
\newtheorem{re}[thm]{Remark}
\newtheorem{exa}[thm]{Example}
\begin{document}
\title {\bf Multiple zeta values and Euler sums}
\author{
{Ce Xu\thanks{Corresponding author. Email: xuce1242063253@163.com (C. XU)}}\\[1mm]
\small School of Mathematical Sciences, Xiamen University\\
\small Xiamen
361005, P.R. China}
\date{}
\maketitle \noindent{\bf Abstract } In this paper, we establish some expressions of series involving harmonic numbers and Stirling numbers of the first kind in terms of multiple zeta values, and present some new relationships between multiple zeta values and multiple zeta star values. The relationships obtained allow us to find some nice closed form representations of nonlinear Euler sums through Riemann zeta values and linear sums. Furthermore, we show that the combined sums
\[H\left( {a,b;m,p} \right) := \sum\limits_{a + b = m - 1} {\zeta \left( {{{\left\{ p \right\}}_a},p + 1,{{\left\{ p \right\}}_b}} \right)}\quad (m\in \N,p>1) \]
and
\[{H^ \star }\left( {a,b;m,p} \right) := \sum\limits_{a + b = m - 1} {{\zeta ^ \star }\left( {{{\left\{ p \right\}}_a},p + 1,{{\left\{ p \right\}}_b}} \right)}\quad (m\in \N,p>1) \]
are reducible to polynomials in zeta values, and give explicit recurrence formulas. Some interesting (known or new) consequences and illustrative examples are considered.
\\[2mm]
\\[2mm]
\section{Introduction}
Let $\R$ and $\mathbb{C}$ denote, respectively the sets of real and complex numbers and let $\N:=\{1,2,3,\ldots\}$ be the set of natural numbers, and $\N_0:=\N\cup \{0\}$ be the set of positive integers and $\mathbb{N} \setminus \{1\}:=\{2,3,4,\cdots\}$. For any multi-index $\mathbf{S} := \left( {{s_1},{s_2}, \cdots ,{s_k}} \right)\ \left( {{s_i} \in \mathbb{N},\ k\in \N,\;{s_1} > 1} \right)$, the general multiple zeta value $\zeta(\mathbf{S})$ and the multiple zeta star value $\zeta^\star(\mathbf{S})$ are defined, respectively, by convergent series (\cite{BBBL1997,Fa2016,KO2010,Y2009})
\[\zeta \left( \mathbf{S} \right) \equiv \zeta \left( {{s_1},{s_2}, \cdots ,{s_k}} \right) := \sum\limits_{{n_1} > {n_2} >  \cdots  > {n_k} \ge 1} {\frac{1}{{n_1^{{s_1}}n_2^{{s_2}} \cdots n_k^{{s_k}}}}} ,\tag{1.1}\]
\[\zeta^\star\left( \mathbf{S} \right) \equiv\zeta^\star \left( {{s_1},{s_2}, \cdots ,{s_k}} \right) := \sum\limits_{{n_1} \ge {n_2} \ge  \cdots  \ge {n_k} \ge 1} {\frac{1}{{n_1^{{s_1}}n_2^{{s_2}} \cdots n_k^{{s_k}}}}}.\tag{1.2}\]
where ${s_1} +  \cdots  + {s_k}$ is called the weight and $k$ is the multiplicity.
For convenience, we let $\{a\}_k$ be the $k$ repetitions of a such that
\[\zeta \left( {5,3,{{\left\{ 1 \right\}}_2}} \right) = \zeta \left( {5,3,1,1} \right),\;{\zeta ^ \star }\left( {4,2,{{\left\{ 1 \right\}}_3}} \right) = {\zeta ^ \star }\left( {4,2,1,1,1} \right).\]
Many papers use the opposite convention, with the $n_i$'s ordered by $n_1<n_2<\cdots<n_k$ or $n_1\leq n_2\leq \cdots\leq n_k$, see \cite{CCE2016,CML2016,Dr1991,E2009,H1992}.
Multiple zeta values and multiple zeta star values were introduced and studied by Euler \cite{E1927} in the old days.  The multiple zeta values have attracted considerable interest in recent years. In the past two decades, many authors have studied multiple zeta values and multiple zeta star values, and a number of relations among them have been found \cite{BBBL1997,BBL2001,BBBL2001,BG1996,CCE2016,CML2016,Dr1991,E2009,E2016,H1992,I2001,L2012,KO2010,CM1994,O1999,KP2013,DZ1994,DZ2012,JZ2010} . There are important properties for multiple zeta values, so called sum, cyclic sum, and duality formulas. For example, one of the well known Q-linear
relations among multiple zeta values is the sum formula ( see\cite{BBBL1997,I2001}), which states that
\[\sum\limits_{\scriptstyle {s_1} +  \cdots+ {s_k} = n \hfill \atop
  \scriptstyle {\rm Each}\;{s_j} \ge 1,\;{s_1} > 1 \hfill} {\zeta \left( {{s_1},{s_{2}}, \cdots ,{s_k}} \right)}  = \zeta \left( n \right).\tag{1.3}\]
In \cite{I2001}, M. Igarashi proved a generalization of the sum formula (see Proposition 3 in the reference \cite{I2001}).
From \cite{BBBL1997,CCE2016,CML2016,Dr1991,E2009,E2016}, we know that multiple zeta values can be represented by iterated integrals
(or Drinfeld integrals) over a simplex of weight dimension. Thus, we have the alternative $(s_1+s_2+\cdots+s_k)$-dimensional iterated-integral representation
\[\zeta \left( {{s_1},{s_{2,}} \cdots ,{s_k}} \right) = \int\limits_0^1 {{\Omega ^{{s_1} - 1}}{w_1}{\Omega ^{{s_2} - 1}}{w_2} \cdots {\Omega ^{{s_k} - 1}}{w_k}} ,\;{s_1} > 1,\tag{1.4}\]
in which the integrand denotes a string of distinct differential 1-forms of type $\Omega  := dx/x,$ and $w_j$ is given by \[{w_j} := \frac{{d{x_j}}}{{1 - {x_j}}}.\tag{1.5}\]
By using (1.4), Jonathan M. Borwein, David M. Bradley and David J. Broadhurst [5] proved the following duality relation
\[\zeta \left( {{m_1} + 2,{{\left\{ 1 \right\}}_{{n_1}}}, \ldots ,{m_p} + 2,{{\left\{ 1 \right\}}_{{n_p}}}} \right) = \zeta \left( {{n_p} + 2,{{\left\{ 1 \right\}}_{{m_p}}}, \ldots ,{n_1} + 2,{{\left\{ 1 \right\}}_{{m_1}}}} \right).\tag{1.6}\]
A generalization of this duality formula can be found in \cite{CCE2016,CML2016,E2009,E2016}.
On the other hand, the corresponding property of the duality formula for multiple zeta-star values was not known until recently. The best result to date are due to Masanobu Kaneko, Yasuo Ohno \cite{KO2010} and Chika Yamazaki \cite{Y2009}. Kaneko and Ohno proved the following property
\[{\left( { - 1} \right)^k}\zeta^\star \left( {k + 1,{{\left\{ 1 \right\}}_n}} \right) - {\left( { - 1} \right)^n}\zeta^\star \left( {n + 1,{{\left\{ 1 \right\}}_k}} \right) \in \mathbf{Q}\left[ {\zeta \left( 2 \right),\zeta \left( 3 \right), \ldots } \right]\tag{1.7}\]
the right-hand side being the algebra over $\mathbf{Q}$ generated by the values of the Riemann zeta function
at positive integer arguments ($s> 1$). C. Yamazaki \cite{Y2009} gave another proof of (1.7). In \cite{KPZ2016}, Hessami
Pilehroods and Zhao found another type of results on general duality
relations between multiple zeta-star values and Euler sums, see Theorem 1.4 in \cite{KPZ2016}.

The subjects of this paper are Multiple zeta values and Euler sums. Next, we give an introduction to the linear and nonlinear Euler sums.
When $k=2$ in (1.2), then $\zeta^\star \left( {{s_1},{s_2}} \right)$ also called the classical linear Euler sums, which is defined by \cite{FS1998}
\[{S_{p,q}} : = \sum\limits_{n = 1}^\infty  {\frac{{{\zeta _n}\left( p \right)}}{{{n^q}}}}= \zeta^\star \left( {q,p} \right) ,\ p\in \N, q\in\mathbb{N} \setminus \{1\},\tag{1.8}\]
where $\zeta_n{(p)}$ stands for the generalized harmonic number defined by
\[{\zeta _n}\left( p \right) := \sum\limits_{j = 1}^n {\frac{1}{{{j^p}}}},p>0,n\in \N,\tag{1.9} \]
when $p=1$, ${H_n} := {\zeta _n}\left( 1 \right)$ is classical harmonic number, the empty sum $\zeta_0{(p)}$ is conventionally understood to be zero (Many papers use the notation $H^{(p)}_n$ to stands for the generalized harmonic number, namely $H^{(p)}_n\equiv \zeta_n{(p)}$). The generalized harmonic number converges to the Riemann zeta function $\zeta(s)$:
 \[\mathop {\lim }\limits_{n \to \infty } \zeta _n{\left( p \right)} = \zeta \left( p \right)\quad (p\in \mathbb{C}, {\Re} \left( p \right) > 1)\]
where the Riemann zeta function is defined by
\[\zeta(p):=\sum\limits_{n = 1}^\infty {\frac {1}{n^{p}}},\Re(p)>1.\tag{1.10}\]
In general, the multiple harmonic number (also called the partial sums of multiple zeta value) and multiple harmonic star number (also called the partial sums of multiple zeta star value) are defined by
\[{\zeta _n}\left( {{s_1},{s_2}, \cdots ,{s_k}} \right): = \sum\limits_{n \ge {n_1} > {n_2} >  \cdots  > {n_k} \ge 1} {\frac{1}{{n_1^{{s_1}}n_2^{{s_2}} \cdots n_k^{{s_k}}}}} ,\tag{1.11}\]
\[{\zeta_n ^ \star }\left( {{s_1},{s_2}, \cdots ,{s_k}} \right): = \sum\limits_{n \ge {n_1} \ge {n_2} \ge  \cdots  \ge {n_k} \ge 1} {\frac{1}{{n_1^{{s_1}}n_2^{{s_2}} \cdots n_k^{{s_k}}}}},\tag{1.12}\]
when $n<k$, then ${\zeta _n}\left( {{s_1},{s_2}, \cdots ,{s_k}} \right)=0$, and ${\zeta _n}\left(\emptyset \right)={\zeta^\star _n}\left(\emptyset \right)=1$.
The generalized (nonlinear) Euler sums are the infinite sums whose general term is a product of harmonic numbers of index $n$ and a power of $n^{-1}$. Namely, for a multi-index ${\bf S}=(s_1,s_2,\ldots,s_k)\ (k,s_i\in \N, i=1,2,\ldots,k)$ with $s_1\leq s_2\leq \ldots\leq s_k$ and $q\geq 2$, the nonlinear Euler sums of index ${\bf S},q$ is defined by (see \cite{FS1998})
\[{S_{{\bf S},q}} := \sum\limits_{n = 1}^\infty  {\frac{\zeta_n(s_1)\zeta_n(s_2)\cdots\zeta_n(s_k)}
{{{n^q}}}},\tag{1.13}\]
where the quantity $w:={s _1} +  \cdots  + {s _k} + q$ is called the weight, the quantity $k$ is called the degree.
As usual, repeated summands in partitions are indicated by powers, so that for instance
\[{S_{{1^2}{2^3}4,q}} = {S_{112224,q}} = \sum\limits_{n = 1}^\infty  {\frac{{H_n^2\zeta _n^3\left( 2 \right){\zeta _n}\left( 4 \right)}}{{{n^q}}}}. \]
It has been discovered in the course of the years that many Euler sums admit expressions involving finitely the ``zeta values", that is to say values of the Riemann zeta function $\zeta(s)$ with the positive integer arguments. Euler started this line of investigation in the course of a correspondence with Goldbach beginning and he was the first to consider the linear sums $S_{p,q}$. Euler showed this problem in the case $p = 1$ and gave a general formula for odd weight $p + q$  in 1775. Moreover, he conjectured that the double linear
sums would be reducible to zeta values when $p + q$ is odd, and even gave what he hoped to obtain the general
formula. In \cite{BBG1995}, D. Borwein, J.M. Borwein and R. Girgensohn proved conjecture and formula, and in \cite{BBG1994}, D.H. Bailey, J.M. Borwein and R. Girgensohn conjectured that the double linear sums when $p + q > 7,p + q$ is even, are not reducible. Hence, the evaluation of $S_{p,q}$ in terms of values of Riemann zeta function at positive integers is known when $p=1,\ p=q,\ (p,q)=(2,4), (4,2)$ or $p+q$ is odd \cite{BBG1994,BBG1995,FS1998}. Similarly, the relationship between the values of the zeta values and nonlinear Euler sums also has been studied by many authors, see \cite{BBG1994,Bl1999,BBG1995,BBGP1996,BB2001,BZB2008,DB2015,FS1998,F2005,J1991,M2010,M2014,R2002,Xu2016,X2016} and references therein. For example, in \cite{FS1998}, Philippe Flajolet and Bruno Salvy proved that all Euler sums of the form $S_{1^p,q}$ for weights $p+q\in \{3,4,5,6,7,9\}$ are expressible polynomially in terms of zeta values and gave explicit formula. For weight 8, all such sums are the sum of a polynomial in zeta values and a rational multiple of $S_{2,6}$, but not the formula. In \cite{X2016,X2017}, we showed that all quadratic Euler sums $S_{1p,q}$ can be evaluated in terms of zeta values and linear sums whenever
$p+q\leq 9,\ p\in\N,\ q\in \mathbb{N} \setminus \{1\}$ and gave explicit formula.

The main purpose of this paper is to establish some relationships between nonlinear Euler sums and multiple zeta values by using the method of iterated integral representations of series. We then use these relations to evaluate several
series with harmonic numbers. Moreover, we prove that all Euler sums of weight $\leq 8$ are reducible to $\mathbb{Q}$-linear combinations of single zeta monomials with the addition of $\{S_{2,6}\}$ for weight 8. Furthermore, we also obtain some explicit formulas between multiple zeta star values and multiple zeta values.
\section{Main results and proofs}
The following lemma will be useful in the development of the main theorems.
\begin{lem}(see \cite{X2016})
Let $m,k$ be integers with $m \ge 2, k \ge 2$, we have the recurrence relation
\[W\left( {m{\rm{,}}k} \right){\rm{ = }} - {\psi ^{\left( {m + k} \right)}}\left( 1 \right)/\left( {k + 1} \right) - \sum\limits_{i = 1}^{m - 1} {\sum\limits_{j = 0}^{k - 1} {\left( {\begin{array}{*{20}{c}}
   {m - 1}  \\
   i  \\
\end{array}} \right)} } \left( {\begin{array}{*{20}{c}}
   k  \\
   j  \\
\end{array}} \right)W\left( {i{\rm{,}}j} \right){\psi ^{\left( {m + k - i - j - 1} \right)}}\left( 1 \right),\tag{2.1}\]
where the integral $W\left( {m{\rm{,}}k} \right)$ is defined by \[W\left( {m{\rm{,}}k} \right){\rm{ := }}\int\limits_0^1 {\frac{{\ln ^k{{\left( {1 - x} \right)}}{{\left( {\ln x} \right)}^m}}}{{1 - x}}} dx, \]
and $\psi \left( z \right)$ stands  for digamma function (or called Psi function) defined by\[\psi \left( z \right):= \frac{d}{{dz}}\left( {\ln \Gamma \left( z \right)} \right) = \frac{{\Gamma '\left( z \right)}}{{\Gamma \left( z \right)}},\]
which is the logarithmic derivative of the Euler gamma function,
when $n\in \N$, then ${\psi ^{\left( n \right)}}\left( z \right) = {\left( { - 1} \right)^{n + 1}}n!\sum\limits_{k = 0}^\infty  {1/{{\left( {z + k} \right)}^{n + 1}}} $, and $\Gamma \left( z \right) := \int\limits_0^\infty  {{e^{ - t}}{t^{z - 1}}dt} ,\;{\mathop{\Re}\nolimits} \left( z \right) > 0$ denotes the gamma function.
By a direct calculation, we can deduce that
\[W\left( {m,0} \right) = {\left( { - 1} \right)^m}m!\zeta \left( {m + 1} \right),W\left( {1,k} \right) = {\left( { - 1} \right)^{k + 1}}k!\zeta \left( {k + 2} \right),\]
\end{lem}
From the recurrence relation (3.1), we know that the values of integrals $W(m,k)$ can be expressed as a rational linear combination of products of zeta values. Here are a few explicit evaluations.
\begin{align*}
&W\left( {2,1} \right) =  - \frac{1}{2}\zeta \left( 4 \right),\\
&W\left( {3,1} \right) = 12\zeta \left( 5 \right) - 6\zeta \left( 2 \right)\zeta \left( 3 \right),\\
&W\left( {4,1} \right) = 12{\zeta ^2}\left( 3 \right) - 18\zeta \left( 6 \right),\\
&W\left( {5,1} \right) = 360\zeta (7) - 120\zeta (3)\zeta (4) - 120\zeta (2)\zeta (5),\\
&W\left( {4,2} \right) = 240\zeta \left( 7 \right) - 60\zeta \left( 3 \right)\zeta \left( 4 \right) - 96\zeta \left( 2 \right)\zeta \left( 5 \right),\\
&W\left( {3,3} \right) = 180\zeta \left( 7 \right) - 45\zeta \left( 3 \right)\zeta \left( 4 \right) - 72\zeta \left( 2 \right)\zeta \left( 5 \right),\\
&W\left( {4,3} \right) =  - \frac{{1497}}{4}\zeta \left( 8 \right) + 576\zeta \left( 3 \right)\zeta \left( 5 \right) - 144\zeta \left( 2 \right){\zeta ^2}\left( 3 \right),\\
&W\left( {3,4} \right) =  - 366\zeta \left( 8 \right) + 432\zeta \left( 3 \right)\zeta \left( 5 \right) - 72\zeta \left( 2 \right){\zeta ^2}\left( 3 \right),\\
&W\left( {5,2} \right) =  - 610\zeta \left( 8 \right) + 720\zeta \left( 3 \right)\zeta \left( 5 \right) - 120\zeta \left( 2 \right){\zeta ^2}\left( 3 \right),\\
&W\left( {6,2} \right) = 13440\zeta \left( 9 \right) + 240{\zeta ^3}\left( 3 \right) - 4320\zeta \left( 2 \right)\zeta \left( 7 \right) - 2520\zeta \left( 3 \right)\zeta \left( 6 \right) - 3240\zeta \left( 4 \right)\zeta \left( 5 \right),\\
&W\left( {4,4} \right) = 8064\zeta \left( 9 \right) + 288{\zeta ^3}\left( 3 \right) - 2880\zeta \left( 2 \right)\zeta \left( 7 \right) - 1260\zeta \left( 3 \right)\zeta \left( 6 \right) - 2016\zeta \left( 4 \right)\zeta \left( 5 \right),\\
&W\left( {3,5} \right) = 6720\zeta \left( 9 \right) + 120{\zeta ^3}\left( 3 \right) - 2160\zeta \left( 2 \right)\zeta \left( 7 \right) - 1260\zeta \left( 3 \right)\zeta \left( 6 \right) - 1620\zeta \left( 4 \right)\zeta \left( 5 \right),\\
&W\left( {5,3} \right) = 10080\zeta \left( 9 \right) + 360{\zeta ^3}\left( 3 \right) - 3600\zeta \left( 2 \right)\zeta \left( 7 \right) - 1575\zeta \left( 3 \right)\zeta \left( 6 \right) - 2520\zeta \left( 4 \right)\zeta \left( 5 \right),\\
&W\left( {5,4} \right) =  - \frac{{84483}}{4}\zeta \left( {10} \right) - 11520\zeta \left( 2 \right)\zeta \left( 3 \right)\zeta \left( 5 \right) + 28800\zeta \left( 3 \right)\zeta \left( 7 \right) - 3600\zeta \left( 4 \right){\zeta ^2}\left( 3 \right) + 14400{\zeta ^2}\left( 5 \right),\\
&W\left( {4,5} \right) =  - 17514\zeta \left( {10} \right) - 8640\zeta \left( 2 \right)\zeta \left( 3 \right)\zeta \left( 5 \right) + 23040\zeta \left( 3 \right)\zeta \left( 7 \right) - 3240\zeta \left( 4 \right){\zeta ^2}\left( 3 \right) + 11520{\zeta ^2}\left( 5 \right).
\end{align*}
\begin{thm} For $n,m\in \N$ and $x\in[-1,1)$. Then the following relation holds:
\begin{align*}
\int\limits_0^x {{t^{n - 1}}{{\ln }^m}\left( {1 - t} \right)} dt&=\frac{1}{n}{\ln ^m}\left( {1 - x} \right)\left( {{x^n} - 1} \right)+ m!\frac{{{{\left( { - 1} \right)}^m}}}{n}\sum\limits_{1 \le {k_m} \le  \cdots  \le {k_1} \le n} {\frac{{{x^{{k_m}}}}}{{{k_1} \cdots {k_m}}}} \\
&\quad - \frac{1}{n}\sum\limits_{i = 1}^{m - 1} {{{\left( { - 1} \right)}^{i - 1}}i!\left( {\begin{array}{*{20}{c}}
   m  \\
   i  \\
\end{array}} \right){{\ln }^{m - i}}\left( {1 - x} \right)} \sum\limits_{1 \le {k_i} \le  \cdots  \le {k_1} \le n} {\frac{{{x^{{k_i}}} - 1}}{{{k_1} \cdots {k_i}}}} .\tag{2.2}
\end{align*}
\end{thm}
\pf The proof is by induction on $m$.
Define $J\left( {n,m;x} \right): = \int\limits_0^x {{t^{n - 1}}{{\ln }^m}\left( {1 - t} \right)} dt$, for $m=1$, by simple calculation, we can arrive at the conclusion that
\[J\left( {n,1;x} \right) = \int_0^x {{t^{n - 1}}\ln \left( {1 - t} \right)} dt = \frac{1}{n}\left\{ {{x^n}\ln \left( {1 - x} \right) - \sum\limits_{j = 1}^n {\frac{{{x^j}}}{j}}  - \ln \left( {1 - x} \right)} \right\},\]
and the formula is true. For $m>1$ we proceed as follows.
By using integration by parts we have the following recurrence relation
\[J\left( {n,m;x} \right) = \frac{1}{n}{\ln ^m}\left( {1 - x} \right)\left( {{x^n} - 1} \right) - \frac{m}{n}\sum\limits_{k = 1}^n {J\left( {k,m - 1;x} \right)}.\tag{2.3} \]
Then by the induction hypothesis, we have
\begin{align*}
J\left( {n,m;x} \right) =& \frac{1}{n}{\ln ^m}\left( {1 - x} \right)\left( {{x^n} - 1} \right) + m!\frac{{{{\left( { - 1} \right)}^m}}}{n}\sum\limits_{j = 1}^n {\sum\limits_{1 \le {k_{m - 1}} \le  \cdots  \le {k_1} \le j} {\frac{{{x^{{k_{m - 1}}}}}}{{{k_1} \cdots {k_{m - 1}}}}} } \\
& + \frac{m}{n}\sum\limits_{i = 1}^{m - 2} {{{\left( { - 1} \right)}^{i - 1}}i!\left( {\begin{array}{*{20}{c}}
   {m - 1}  \\
   i  \\
\end{array}} \right){{\ln }^{m - 1 - i}}\left( {1 - x} \right)} \sum\limits_{1 \le {k_i} \le  \cdots  \le {k_1} \le j} {\frac{{{x^{{k_i}}} - 1}}{{{k_1} \cdots {k_i}}}} \\
& - \frac{m}{n}{\ln ^{m - 1}}\left( {1 - x} \right)\sum\limits_{j = 1}^n {\frac{{{x^j} - 1}}{j}} .\tag{2.4}
\end{align*}
Thus, by a direct calculation, we deduce the desired result. This completes the proof of Theorem 2.2. \hfill$\square$\\
Letting $x$ approach 1 in (2.2) and using the definition of multiple star harmonic number, we obtain
\begin{align*}
&\int_0^1 {{t^{n - 1}}\ln^m \left( {1 - t} \right)} dt = m!\frac{{{{\left( { - 1} \right)}^m}}}{n}{\zeta^\star _n}\left( {{{\left\{ 1 \right\}}_m}} \right),\tag{2.5}
\end{align*}
Dividing (2.2) by $x$ and integrating over the interval (0,1), the result is
\begin{align*}
 \int\limits_0^1 {{x^{n - 1}}\ln x{{\ln }^m}\left( {1 - x} \right)} dx =& m!\frac{{{{\left( { - 1} \right)}^{m - 1}}}}{n}\left\{ {\frac{{\zeta^\star_n\left( {{{\left\{ 1 \right\}}_m}} \right)}}{n} - \zeta \left( {m + 1} \right)} \right\} + m!\frac{{{{\left( { - 1} \right)}^{m - 1}}}}{n}\zeta^\star_n\left( {{{\left\{ 1 \right\}}_{m - 1}},2} \right) \\
  &+ m!\frac{{{{\left( { - 1} \right)}^{m - 1}}}}{n}\sum\limits_{i = 1}^{m - 1} {\left\{ {\zeta^\star _n\left( {{{\left\{ 1 \right\}}_{i - 1}},2,{{\left\{ 1 \right\}}_{m - i}}} \right) - \zeta \left( {m - i + 1} \right)\zeta^\star _n\left( {{{\left\{ 1 \right\}}_i}} \right)} \right\}}.\tag{2.6}
\end{align*}
On the other hand, in \cite{X2016}, we proved the identity
\[\int\limits_0^1 {{t^{n - 1}}{{\ln }^k}\left( {1 - t} \right)} dt = {\left( { - 1} \right)^k}\frac{{{Y_k}\left( n \right)}}{n},\ {Y_0}\left( n \right) = 1,\tag{2.7}\]
where ${Y_k}\left( n \right) := {Y_k}\left( {{\zeta _n}\left( 1 \right),1!{\zeta _n}\left( 2 \right),2!{\zeta _n}\left( 3 \right), \cdots ,\left( {r - 1} \right)!{\zeta _n}\left( r \right), \cdots } \right)$, and ${Y_k}\left( {{x_1},{x_2}, \cdots } \right)$ stands for the complete exponential Bell polynomial defined by (see \cite{L1974})
\[\exp \left( {\sum\limits_{m \ge 1}^{} {{x_m}\frac{{{t^m}}}{{m!}}} } \right) = 1 + \sum\limits_{k \ge 1}^{} {{Y_k}\left( {{x_1},{x_2}, \cdots } \right)\frac{{{t^k}}}{{k!}}}.\tag{2.8}\]
The relations (2.5) and (2.7) yield the following result
\[\zeta _n^ \star \left( {{{\{ 1\} }_m}} \right) = \frac{1}{{m!}}{Y_m}\left( n \right),n,m\in \N_0.\tag{2.9}\]
By using the definition of the complete exponential Bell polynomial, we easily deduce that
\begin{align*}
&{Y_1}\left( n \right) = {H_n},\\&{Y_2}\left( n \right) = H_n^2 + {\zeta _n}\left( 2 \right),\\&{Y_3}\left( n \right) =  H_n^3+ 3{H_n}{\zeta _n}\left( 2 \right)+ 2{\zeta _n}\left( 3 \right),\\
&{Y_4}\left( n \right) = H_n^4 + 8{H_n}{\zeta _n}\left( 3 \right) + 6H_n^2{\zeta _n}\left( 2 \right) + 3\zeta _n^2\left( 2 \right) + 6{\zeta _n}\left( 4 \right),\\
&{Y_5}\left( n \right) = \left\{ \begin{array}{l}
 H_n^5 + 10H_n^3{\zeta _n}\left( 2 \right) + 20H_n^2{\zeta _n}\left( 3 \right) + 15{H_n}\zeta _n^2\left( 2 \right) \\
  + 30{H_n}{\zeta _n}\left( 4 \right) + 20{\zeta _n}\left( 2 \right){\zeta _n}\left( 3 \right) + 24{\zeta _n}\left( 5 \right) \\
 \end{array} \right\},\\
&{Y_6}\left( n \right) = \left\{ \begin{array}{l}
 H_n^6 + 15H_n^4{\zeta _n}\left( 2 \right) + 40H_n^3{\zeta _n}\left( 3 \right) + 90H_n^2{\zeta _n}\left( 4 \right) \\
  + 144{H_n}{\zeta _n}\left( 5 \right) + 45H_n^2\zeta _n^2\left( 2 \right) + 120{H_n}{\zeta _n}\left( 2 \right){\zeta _n}\left( 3 \right) \\
  + 40\zeta _n^2\left( 3 \right) + 15\zeta _n^3\left( 2 \right) + 90{\zeta _n}\left( 2 \right){\zeta _n}\left( 4 \right) + 120{\zeta _n}\left( 6 \right) \\
 \end{array} \right\},
\end{align*}
and we can obtain the conclusion: For any $n,m\in \N$, ${Y_k}\left( n \right)$ is a rational linear combination of products of harmonic numbers. Furthermore, we know that the multiple star harmonic number $\zeta _n^ \star \left( {{{\left\{ 1 \right\}}_m}} \right)$ can be expressed as a rational linear combination of harmonic numbers. In fact, using Eulerian beta integral, we can show the more general recurrence relation:
\begin{thm}
For integers $m,k,n\in \N$, then we have the following recurrence relation
\begin{align*}
&I\left( {n,m,k} \right) = \sum\limits_{i = 0}^{m - 1} {\left( {\begin{array}{*{20}{c}}
   {m - 1}  \\
   i  \\
\end{array}} \right)\left( {m - i - 1} \right)!\frac{{{{\left( { - 1} \right)}^{m - i}}}}{{{n^{m - i}}}}I\left( {n,i,k} \right)} \\
&\quad\  + \sum\limits_{j = 0}^{k - 1} {\left( {\begin{array}{*{20}{c}}
   k  \\
   j  \\
\end{array}} \right){{\left( { - 1} \right)}^{m + k - j}}\left( {m + k - j - 1} \right)!{\zeta _n}\left( {m + k - j} \right)I\left( {n,0,j} \right)}  \\
&\quad\  - \sum\limits_{j = 0}^{k - 1} {\left( {\begin{array}{*{20}{c}}
   k  \\
   j  \\
\end{array}} \right){{\left( { - 1} \right)}^{m + k - j}}\left( {m + k - j - 1} \right)!\zeta \left( {m + k - j} \right)I\left( {n,0,j} \right)}\\
&\quad\  + \sum\limits_{i = 1}^{m - 1} {\sum\limits_{j = 0}^{k - 1} {\left( {\begin{array}{*{20}{c}}
   {m - 1}  \\
   i  \\
\end{array}} \right)\left( {\begin{array}{*{20}{c}}
   k  \\
   j  \\
\end{array}} \right){{\left( { - 1} \right)}^{m + k - i - j}}\left( {m + k - i - j - 1} \right)!{\zeta _n}\left( {m + k - i - j} \right)I\left( {n,i,j} \right)} } \\
&\quad\  - \sum\limits_{i = 1}^{m - 1} {\sum\limits_{j = 0}^{k - 1} {\left( {\begin{array}{*{20}{c}}
   {m - 1}  \\
   i  \\
\end{array}} \right)\left( {\begin{array}{*{20}{c}}
   k  \\
   j  \\
\end{array}} \right){{\left( { - 1} \right)}^{m + k - i - j}}\left( {m + k - i - j - 1} \right)!\zeta \left( {m + k - i - j} \right)I\left( {n,i,j} \right)} }.\tag{2.10}
\end{align*}
where $I\left( {n ,m,k} \right)$ is defined by the integral
\[I\left( {n ,m,k} \right): = \int\limits_0^1 {{x^{n  - 1}}{{\ln }^m}x{{\ln }^k}\left( {1 - x} \right)} dx.\tag{2.11}\]
\[I\left( {n ,0,0} \right) = \frac{1}{n },I\left( {n ,i,0} \right) = {\left( { - 1} \right)^i}i!\frac{1}{{{n ^{i + 1}}}}.\]
When $m=0$, then $I\left( {n,0,k} \right) = J\left( {n,k,1} \right)$.
\end{thm}
\pf Applying the definition of Beta function ${B\left( {\alpha ,\beta } \right)}$, we can find that
\[I\left( {n ,m,k} \right): = \int\limits_0^1 {{x^{n  - 1}}{{\ln }^m}x{{\ln }^k}\left( {1 - x} \right)} dx = {\left. {\frac{{{\partial ^{m + k}}B\left( {\alpha ,\beta } \right)}}{{\partial {\alpha ^m}\partial {\beta ^k}}}} \right|_{\alpha=n,\beta  = 1}},\tag{2.12}\]
where the Eulerian Beta function is defined by
\[B\left( {\alpha,\beta} \right) := \int\limits_0^1 {{x^{\alpha - 1}}{{\left( {1 - x} \right)}^{\beta - 1}}dx}  = \frac{{\Gamma \left( \alpha \right)\Gamma \left( \beta \right)}}{{\Gamma \left( {\alpha + \beta} \right)}},\;{\mathop{\Re}\nolimits} \left( \alpha \right) > 0,{\mathop{\Re}\nolimits} \left( \beta \right) > 0.\tag{2.13}\]
By using (2.13) and the definition of digamma function $\psi (x)$, it is obvious that
\[\frac{{\partial B\left( {\alpha ,\beta } \right)}}{{\partial \alpha }} = B\left( {\alpha ,\beta } \right)\left[ {\psi \left( \alpha  \right) - \psi \left( {\alpha  + \beta } \right)} \right].\]
Therefore, differentiating $m-1$ times this equality, we can deduce that
\[\frac{{{\partial ^m}B\left( {\alpha ,\beta } \right)}}{{\partial {\alpha ^m}}} = \sum\limits_{i = 0}^{m - 1} {\left( {\begin{array}{*{20}{c}}
   {m - 1}  \\
   i  \\
\end{array}} \right)\frac{{{\partial ^i}B\left( {\alpha ,\beta } \right)}}{{\partial {\alpha ^i}}}} \cdot\left[ {{\psi ^{\left( {m - i - 1} \right)}}\left( \alpha  \right) - {\psi ^{\left( {m - i - 1} \right)}}\left( {\alpha  + \beta } \right)} \right].\tag{2.14}\]
Since $B(\alpha,\beta)=B(\beta,\alpha)$, the change of variable $\alpha\mapsto \beta,\ \beta\mapsto\alpha$, then we also have
\[\frac{{{\partial ^m}B\left( {\alpha ,\beta } \right)}}{{\partial {\beta ^m}}} = \sum\limits_{i = 0}^{m - 1} {\left( {\begin{array}{*{20}{c}}
   {m - 1}  \\
   i  \\
\end{array}} \right)\frac{{{\partial ^i}B\left( {\alpha ,\beta } \right)}}{{\partial {\beta ^i}}}} \cdot\left[ {{\psi ^{\left( {m - i - 1} \right)}}\left( \beta  \right) - {\psi ^{\left( {m - i - 1} \right)}}\left( {\beta  + \alpha } \right)} \right].\tag{2.15}\]
Putting $\alpha=n, \beta=1$ in (2.15), we arrive at the conclusion that
\[I\left( {n ,0,m} \right){\rm{ = }}\sum\limits_{i = 0}^{m - 1} {{{\left( { - 1} \right)}^{m - i}}\left( {m - i - 1} \right)!\left( {\begin{array}{*{20}{c}}
   {m - 1}  \\
   i  \\
\end{array}} \right)I\left( {n ,0,i} \right)} \zeta_n {\left( {m - i} \right)}.\tag{2.16}\]
Furthermore, by using (2.14), the following identity is easily derived
\begin{align*}
\frac{{{\partial ^{m + k}}B\left( {\alpha,\beta} \right)}}{{\partial {\alpha^m}\partial {\beta^k}}}
 &=\frac{{{\partial ^k}}}{{\partial {\beta^k}}}\left( {\frac{{{\partial ^m}B\left( {\alpha,\beta} \right)}}{{\partial {\alpha^m}}}} \right)
\nonumber \\
           &=\sum\limits_{i = 0}^{m - 1} {\left( {\begin{array}{*{20}{c}}
   {m - 1}  \\
   i  \\
\end{array}} \right)\frac{{{\partial ^{i + k}}B\left( {\alpha,\beta} \right)}}{{\partial {\alpha^i}\partial {\beta^k}}} \cdot } \left[ {{\psi ^{\left( {m - i - 1} \right)}}\left(\alpha \right) - {\psi ^{\left( {m - i - 1} \right)}}\left( {\alpha + \beta} \right)} \right]
\nonumber \\
           &\quad \  - \sum\limits_{j = 0}^{k - 1} {\left( {\begin{array}{*{20}{c}}
   k  \\
   j  \\
\end{array}} \right)} \frac{{{\partial ^{ j}}B\left( {\alpha,\beta} \right)}}{{\partial {\beta^j}}}{\psi ^{\left( {m + k - j - 1} \right)}}\left( {\alpha + \beta} \right)
\nonumber \\
           &\quad \  - \sum\limits_{i = 1}^{m - 1} {\sum\limits_{j = 0}^{k - 1} {\left( {\begin{array}{*{20}{c}}
   {m - 1}  \\
   i  \\
\end{array}} \right)\left( {\begin{array}{*{20}{c}}
   k  \\
   j  \\
\end{array}} \right)} \frac{{{\partial ^{i + j}}B\left( {\alpha,\beta} \right)}}{{\partial {\alpha^i}\partial {\beta^j}}}} {\psi ^{\left( {m + k - i - j - 1} \right)}}\left( {\alpha + \beta} \right)
.\tag{2.17}
\end{align*}
By the definition of digamma function again, we know that
\[{\psi ^{\left( {m - i - 1} \right)}}\left( n  \right) - {\psi ^{\left( {m - i - 1} \right)}}\left( {n  + 1} \right) = {\left( { - 1} \right)^{m - i}}\left( {m - i - 1} \right)!\frac{1}{{{n ^{m - i}}}},\tag{2.18}\]
\[{\psi ^{\left( {m + k - j - 1} \right)}}\left( {n  + 1} \right) = {\left( { - 1} \right)^{m + k - j}}\left( {m + k - j - 1} \right)!\left( {\zeta \left( {m + k - j} \right) - \zeta_n {\left( {m + k - j} \right)}} \right).\tag{2.19}\]
Hence, taking $\alpha=n,\beta=1$ in (2.17), then substituting (2.18) and (2.19) into (2.17) respectively, we can obtain (2.1). The proof of Theorem 2.3 is finished.\hfill$\square$\\
From (2.10) and (2.16), we give the following identities: for $n\in \N$,
\begin{align*}
&I\left( {n,0,1} \right) = \int\limits_0^1 {{x^{n - 1}}\ln \left( {1 - x} \right)dx}  =  - \frac{{{H_n}}}{n},\\
&I\left( {n,0,2} \right) = \int\limits_0^1 {{x^{n - 1}}{{\ln }^2}\left( {1 - x} \right)dx}  = \frac{{H_n^2 + {\zeta _n}\left( 2 \right)}}{n},\\
&I\left( {n,1,1} \right) = \int\limits_0^1 {{x^{n - 1}}\ln x\ln \left( {1 - x} \right)dx}  = \frac{{{H_n}}}{{{n^2}}} - \frac{{\zeta \left( 2 \right) - {\zeta _n}\left( 2 \right)}}{n},\\
&I\left( {n,0,3} \right) = \int\limits_0^1 {{x^{n - 1}}{{\ln }^3}\left( {1 - x} \right)dx}  =  - \frac{{H_n^3 + 3{H_n}{\zeta _n}\left( 2 \right) + 2{\zeta _n}\left( 3 \right)}}{n},\\
&I\left( {n,1,2} \right) = \int\limits_0^1 {{x^{n - 1}}\ln x{{\ln }^2}\left( {1 - x} \right)dx}  =  - \frac{{H_n^2 + {\zeta _n}\left( 2 \right)}}{{{n^2}}} + 2\frac{{\zeta \left( 3 \right) - {\zeta _n}\left( 3 \right)}}{n} + 2\frac{{\zeta \left( 2 \right) - {\zeta _n}\left( 2 \right)}}{n}{H_n}.
\end{align*}
By using integration by parts, we can find the relation
\[I\left( {n,m,1} \right) =  - \frac{m}{n}I\left( {n,m - 1,1} \right) + {\left( { - 1} \right)^m}m!\frac{{\zeta \left( {m + 1} \right) - {\zeta _n}\left( {m + 1} \right)}}{n}.\tag{2.20}\]
Thus, by a direct calculation, we obtain the result
\[I\left( {n,m,1} \right) = \int\limits_0^1 {{x^{n - 1}}{{\ln }^m}x\ln \left( {1 - x} \right)dx}  = {\left( { - 1} \right)^{m + 1}}m!\frac{{{H_n}}}{{{n^{m + 1}}}} + {\left( { - 1} \right)^m}m!\sum\limits_{j = 1}^m {\frac{{\zeta \left( {j + 1} \right) - {\zeta _n}\left( {j + 1} \right)}}{{{n^{m - j + 1}}}}}.\tag{2.21} \]
Therefore, from Theorem 2.3, we know that the integral $I\left( {n,m,k} \right)$ is a rational linear combination of products of harmonic numbers and zeta values. Taking $m=1$ in (2.10), we can give the the following corollary.
\begin{cor} For integers $n,k\in \N$, then we have
\[I\left( {n,1,k} \right) =  - \frac{1}{n}I\left( {n,0,k} \right) - \sum\limits_{j = 0}^{k - 1} {\left( {\begin{array}{*{20}{c}}
   k  \\
   j  \\
\end{array}} \right){{\left( { - 1} \right)}^{k + 1 - j}}\left( {k - j} \right)!\left( {\zeta \left( {k + 1 - j} \right) - {\zeta _n}\left( {k + 1 - j} \right)} \right)} I\left( {n,0,j} \right).\tag{2.22}\]
\end{cor}
Comparing (2.6) and (2.22), we obtain the conclusion: if $m$ and $n$ are positive integers, then the sums of multiple star harmonic number $\sum\limits_{i = 1}^m {{\zeta_n^ \star }\left( {{{\left\{ 1 \right\}}_{i - 1}},2,{{\left\{ 1 \right\}}_{m - i}}} \right)} $ can be expressed in terms of polynomials of the harmonic numbers, and we have the formula
\[\sum\limits_{i = 1}^m {{\zeta_n ^ \star }\left( {{{\left\{ 1 \right\}}_{i - 1}},2,{{\left\{ 1 \right\}}_{m - i}}} \right)}  = \sum\limits_{i = 1}^m {{\zeta _n}\left( {m + 2 - i} \right){\zeta_n^ \star }\left( {{{\left\{ 1 \right\}}_{i - 1}}} \right)}.\tag{2.23}\]
\begin{thm}
For integer $k>0$ and $x\in [-1,1)$, then we have
\[{\ln ^k}\left( {1 - x} \right) = {\left( { - 1} \right)^k}k!\sum\limits_{n = 1}^\infty  {\frac{{{x^n}}}{n}{\zeta _{n - 1}}\left( {{{\left\{ 1 \right\}}_{k - 1}}} \right)},\tag{2.24} \]
\[s\left( {n,k} \right) = \left( {n - 1} \right)!{\zeta _{n - 1}}\left( {{{\left\{ 1 \right\}}_{k - 1}}} \right).\tag{2.25}\]
where ${s\left( {n,k} \right)}$ denotes the (unsigned) Stirling number of the first kind (see \cite{L1974}).
\begin{align*}
& s\left( {n,1} \right) = \left( {n - 1} \right)!,\\&s\left( {n,2} \right) = \left( {n - 1} \right)!{H_{n - 1}},\\&s\left( {n,3} \right) = \frac{{\left( {n - 1} \right)!}}{2}\left[ {H_{n - 1}^2 - {\zeta _{n - 1}}\left( 2 \right)} \right],\\
&s\left( {n,4} \right) = \frac{{\left( {n - 1} \right)!}}{6}\left[ {H_{n - 1}^3 - 3{H_{n - 1}}{\zeta _{n - 1}}\left( 2 \right) + 2{\zeta _{n - 1}}\left( 3 \right)} \right], \\
&s\left( {n,5} \right) = \frac{{\left( {n - 1} \right)!}}{{24}}\left[ {H_{n - 1}^4 - 6{\zeta _{n - 1}}\left( 4 \right) - 6H_{n - 1}^2{\zeta _{n - 1}}\left( 2 \right) + 3\zeta _{n - 1}^2\left( 2 \right) + 8H_{n - 1}^{}{\zeta _{n - 1}}\left( 3 \right)} \right],\\
&s\left( {n,6} \right) = \frac{{\left( {n - 1} \right)!}}{{120}}\left\{ \begin{array}{l}
 H_{n - 1}^5 - 10H_{n - 1}^3{\zeta _{n - 1}}\left( 2 \right) + 20H_{n - 1}^2{\zeta _{n - 1}}\left( 3 \right) + 15{H_{n - 1}}\zeta _{n - 1}^2\left( 2 \right) \\
  - 30{H_{n - 1}}{\zeta _{n - 1}}\left( 4 \right) - 20{\zeta _{n - 1}}\left( 2 \right){\zeta _{n - 1}}\left( 3 \right) + 24{\zeta _{n - 1}}\left( 5 \right) \\
 \end{array} \right\},\\
&s\left( {n,7} \right) = \frac{{\left( {n - 1} \right)!}}{{720}}\left\{ \begin{array}{l}
 H_{n - 1}^6 - 15H_{n - 1}^4{\zeta _{n - 1}}\left( 2 \right) + 40H_{n - 1}^3{\zeta _{n - 1}}\left( 3 \right) - 90H_{n - 1}^2{\zeta _{n - 1}}\left( 4 \right) \\
  + 144{H_{n - 1}}{\zeta _{n - 1}}\left( 5 \right) + 45H_{n - 1}^2\zeta _{n - 1}^2\left( 2 \right) - 120{H_{n - 1}}{\zeta _{n - 1}}\left( 2 \right){\zeta _{n - 1}}\left( 3 \right) \\
  + 40\zeta _{n - 1}^2\left( 3 \right) - 15\zeta _{n - 1}^3\left( 2 \right) + 90{\zeta _{n - 1}}\left( 2 \right){\zeta _{n - 1}}\left( 4 \right) - 120{\zeta _{n - 1}}\left( 6 \right) \\
 \end{array} \right\}.
\end{align*}
The Stirling numbers ${s\left( {n,k} \right)}$ of the first kind satisfy a recurrence relation in the form
\[s\left( {n,k} \right) = s\left( {n - 1,k - 1} \right) + \left( {n - 1} \right)s\left( {n - 1,k} \right),\;\;n,k \in \N,\]
with $s\left( {n,k} \right) = 0,n < k,s\left( {n,0} \right) = s\left( {0,k} \right) = 0,s\left( {0,0} \right) = 1$.
\end{thm}
\pf To prove the first identity we proceed by induction on $k$. Obviously, it is valid for $k=1$. For $k>1$ using the equality
\[{\ln ^{k{\rm{ + }}1}}\left( {1 - x} \right){\rm{ = }} - \left( {k + 1} \right)\int\limits_0^x {\frac{{{{\ln }^k}\left( {1 - t} \right)}}{{1 - t}}dt} \]
and applying the induction hypothesis, by using the Cauchy product of power series, we arrive at
\begin{align*}
{\ln ^{k{\rm{ + }}1}}\left( {1 - x} \right){\rm{ = }}& - \left( {k + 1} \right)\int\limits_0^x {\frac{{{{\ln }^k}\left( {1 - t} \right)}}{{1 - t}}dt} \\
& = {\left( { - 1} \right)^{k + 1}}\left( {k + 1} \right)!\sum\limits_{n = 1}^\infty  {\frac{1}{{n + 1}}\sum\limits_{i = 1}^n {\frac{{{\zeta _{i - 1}}\left( {{{\left\{ 1 \right\}}_{k - 1}}} \right)}}{i}} } {x^{n + 1}}\\
& = {\left( { - 1} \right)^{k + 1}}\left( {k + 1} \right)!\sum\limits_{n = 1}^\infty  {\frac{{{\zeta _n}\left( {{{\left\{ 1 \right\}}_k}} \right)}}{{n + 1}}} {x^{n + 1}}.
\end{align*}
Nothing that ${\zeta _n}\left( {{{\left\{ 1 \right\}}_k}} \right) = 0$ when $n<k$. Thus,
we can deduce (2.24). To prove the second identity of our theorem, we use the following equation (\cite{L1974})
\[{\ln ^k}\left( {1 - x} \right) = {\left( { - 1} \right)^k}k!\sum\limits_{n = k}^\infty  {\frac{{s\left( {n,k} \right)}}{{n!}}{x^n}} ,\: - 1 \le x < 1.\tag{2.26}\]
Comparing the coefficients of $x^n$ in (2.24) and (2.26), we obtain formula (2.25). The proof of Theorem 2.5 is thus completed.\hfill$\square$\\
From (2.1) and (2.26), we have the result
\begin{align*}
W\left( {k,m} \right)& =\int\limits_0^1 {\frac{{{{\ln }^k}\left( {1 - x} \right){{\ln }^m}x}}{x}} dx
 = {\left( { - 1} \right)^k}k!\sum\limits_{n = 1}^\infty  {\frac{{{\zeta _{n - 1}}\left( {{{\left\{ 1 \right\}}_{k - 1}}} \right)}}{n}\int\limits_0^1 {{x^n}{{\ln }^m}x} } dx\\
& = {\left( { - 1} \right)^{m + k}}m!k!\sum\limits_{n = 1}^\infty  {\frac{{{\zeta _{n - 1}}\left( {{{\left\{ 1 \right\}}_{k - 1}}} \right)}}{{{n^{m + 2}}}}}\\
&={\left( { - 1} \right)^{m + k}}m!k!\zeta \left( {m + 2,{{\left\{ 1 \right\}}_{k - 1}}} \right),\tag{2.27}
\end{align*}
which combined with (2.1), implies that for any $m,k\in \N$, the multiple zeta value $\zeta(m+1,\{1\}_{k-1} )$ can be
represented as a polynomial of zeta values with rational coefficients. For example:
\[\begin{array}{l}
 \zeta \left( {2,{{\left\{ 1 \right\}}_m}} \right) = \zeta \left( {m + 2} \right), \\
 \zeta \left( {3,{{\left\{ 1 \right\}}_m}} \right) = \frac{{m + 2}}{2}\zeta \left( {m + 3} \right) - \frac{1}{2}\sum\limits_{k = 1}^m {\zeta \left( {k + 1} \right)\zeta \left( {m + 2 - k} \right)} . \\
 \end{array}\]
 By integration by parts, we deduce that
 \[W\left( {m,k - 1} \right) = \frac{m}{k}W\left( {k,m - 1} \right).\tag{2.28}\]
Applying (2.27) into (2.28), we obtain the well known duality formula
\[\zeta \left( {k + 1,{{\left\{ 1 \right\}}_{m - 1}}} \right) = \zeta \left( {m + 1,{{\left\{ 1 \right\}}_{k - 1}}} \right).\]
Obviously, this result is a special case of the duality of multiple zeta values (1.6). From (2.25)-(2.28), we give the result
\[\sum\limits_{n = p}^\infty  {\frac{{s\left( {n,k} \right)}}{{n!{n^m}}}}  = \zeta \left( {k + 1,{{\left\{ 1 \right\}}_{m - 1}}} \right) = \zeta \left( {m + 1,{{\left\{ 1 \right\}}_{k - 1}}} \right).\tag{2.29}\]
Now we state our main results. The main results of this section are the following theorems.
\begin{thm} For integers $p\in \N$ and $m\in \N_0$. Then
 \[\sum\limits_{n = 1}^\infty  {\frac{{{H_n}s\left( {n,p} \right)}}{{n!{n^{m + 1}}}}}  = \left( {p + 1} \right)\zeta \left( {p + 2,{{\left\{ 1 \right\}}_m}} \right) + \sum\limits_{i = 1}^m {\zeta \left( {p + 1,{{\left\{ 1 \right\}}_{i - 1}},2,{{\left\{ 1 \right\}}_{m - i}}} \right)}.\tag{2.30} \]
\end{thm}
\pf Applying (2.26), we can show that
\begin{align*}
{\left( { - 1} \right)^{p + 1}}p!\sum\limits_{n = 1}^\infty  {\frac{{{H_n}s\left( {n,p} \right)}}{{n!{n^{m + 1}}}}}  =& \int\limits_0^1 {\frac{{\ln \left( {1 - {t_1}} \right)}}{{{t_1}}}} d{t_1}\int\limits_0^{{t_1}} {\frac{1}{{{t_2}}}} d{t_2} \cdots \int\limits_0^{{t_{m - 1}}} {\frac{1}{{{t_m}}}} d{t_m}\int\limits_0^{{t_m}} {\frac{{{{\ln }^p}\left( {1 - {t_{m + 1}}} \right)}}{{{t_{m + 1}}}}} d{t_{m + 1}}\\
 = &\int\limits_{0 < {t_{m + 1}} < {t_m} <  \cdots  < {t_1} < 1} {\frac{{\ln \left( {1 - {t_1}} \right){{\ln }^p}\left( {1 - {t_{m + 1}}} \right)}}{{{t_1}{t_2} \cdots {t_{m + 1}}}}} d{t_1}d{t_2} \cdots d{t_{m + 1}}.\tag{2.31}
\end{align*}
Applying the change of variables
${t_i} \mapsto 1 - {t_{m + 2 - i}}\quad (i = 1,2, \cdots ,m + 1)$
to the above multiple integral, we get the identity
\begin{align*}
&\int\limits_{0 < {t_{m + 1}} < {t_m} <  \cdots  < {t_1} < 1} {\frac{{\ln \left( {1 - {t_1}} \right){{\ln }^p}\left( {1 - {t_{m + 1}}} \right)}}{{{t_1}{t_2} \cdots {t_{m + 1}}}}} d{t_1}d{t_2} \cdots d{t_{m + 1}}\\
& = \int\limits_{0 < {t_{m + 1}} < {t_m} <  \cdots  < {t_1} < 1} {\frac{{\ln \left( {{t_{m + 1}}} \right){{\ln }^p}\left( {{t_1}} \right)}}{{\left( {1 - {t_{m + 1}}} \right)\left( {1 - {t_m}} \right) \cdots \left( {1 - {t_1}} \right)}}} d{t_1}d{t_2} \cdots d{t_{m + 1}}\\
& = \int\limits_0^1 {\frac{{{{\ln }^p}\left( {{t_1}} \right)}}{{1 - {t_1}}}} d{t_1}\int\limits_0^{{t_1}} {\frac{1}{{1 - {t_2}}}} d{t_2} \cdots \int\limits_0^{{t_{m - 1}}} {\frac{1}{{1 - {t_m}}}} d{t_m}\int\limits_0^{{t_m}} {\frac{{\ln \left( {{t_{m + 1}}} \right)}}{{1 - {t_{m + 1}}}}} d{t_{m + 1}},\\
& = \sum\limits_{{n_{m + 1}} = 1}^\infty  {\sum\limits_{{n_m} = 1}^\infty  { \cdots \sum\limits_{{n_1} = 1}^\infty  {\int\limits_0^1 {t_1^{{n_{m + 1}} - 1}{{\ln }^p}\left( {{t_1}} \right)d{t_1}\int\limits_0^{{t_1}} {t_2^{{n_m} - 1}} d{t_2} \cdots \int\limits_0^{{t_{m - 1}}} {t_m^{{n_2-1}}} d{t_m}} } } } \int\limits_0^{{t_m}} {t_{m + 1}^{{n_1} - 1}\ln \left( {{t_{m + 1}}} \right)} d{t_{m + 1}}.\tag{2.32}
\end{align*}
By using integration by parts, we deduce that, for $n,m\in \N$,
\[\int\limits_0^x {{t^{n - 1}}{{\left( {\ln t} \right)}^m}} dt = \sum\limits_{l = 0}^m {l!\left( {\begin{array}{*{20}{c}}
   m  \\
   l  \\
\end{array}} \right)\frac{{{{\left( { - 1} \right)}^l}}}{{{n^{l + 1}}}}{{\left( {\ln x} \right)}^{m - l}}{x^n}},x\in (0,1),\tag{2.33} \]
Taking $m=2$ in (2.33), we get
\[\int\limits_0^x {{t^{n - 1}}\ln t} dt = \frac{1}{n}{x^n}\ln x - \frac{{{x^n}}}{{{n^2}}}.\tag{2.34}\]
Substituting (2.34) into (2.32) and combining (2.31), we arrive at the conclusion that
\begin{align*}
&{\left( { - 1} \right)^{p + 1}}p!\sum\limits_{n = 1}^\infty  {\frac{{{H_n}s\left( {n,p} \right)}}{{n!{n^{m + 1}}}}}\\  = &\sum\limits_{{n_1},{n_2}, \cdots {n_{m + 1}} = 1}^\infty  {\frac{{\left( {p + 1} \right)!{{\left( { - 1} \right)}^{p + 1}}}}{{{n_1}\left( {{n_1} + {n_2}} \right) \cdots \left( {{n_1} +  \cdots  + {n_m}} \right){{\left( {{n_1} +  \cdots  + {n_m} + {n_{m + 1}}} \right)}^{p + 2}}}}} \\
& - \sum\limits_{{n_1},{n_2}, \cdots {n_{m + 1}} = 1}^\infty  {\frac{{p!{{\left( { - 1} \right)}^p}}}{{{n_1}\left( {{n_1} + {n_2}} \right) \cdots {{\left( {{n_1} +  \cdots  + {n_m}} \right)}^2}{{\left( {{n_1} +  \cdots  + {n_m} + {n_{m + 1}}} \right)}^{p + 1}}}}} \\
& -  \cdots \\
&- \sum\limits_{{n_1},{n_2}, \cdots {n_{m + 1}} = 1}^\infty  {\frac{{{{\left( { - 1} \right)}^p}p!}}{{n_1^2\left( {{n_1} + {n_2}} \right) \cdots \left( {{n_1} +  \cdots  + {n_m}} \right){{\left( {{n_1} +  \cdots  + {n_m} + {n_{m + 1}}} \right)}^{p + 1}}}}} \\
 =& {\left( { - 1} \right)^{p + 1}}\left( {p + 1} \right)!\zeta \left( {p + 2,{{\left\{ 1 \right\}}_m}} \right) + {\left( { - 1} \right)^{p + 1}}p!\sum\limits_{i = 1}^m {\zeta \left( {p + 1,{{\left\{ 1 \right\}}_{i - 1}},2,{{\left\{ 1 \right\}}_{m - i}}} \right)}.\tag{2.35}
\end{align*}
Thus, this completes the proof of Theorem 2.6.\hfill$\square$
\begin{thm} For integers $m,p\in \N$ and $r\in \mathbb{N} \setminus \{1\}$. Then
\[\sum\limits_{n = 1}^\infty  {\frac{{s\left( {n,m} \right){\zeta _n}\left( r \right)}}{{n!{n^p}}}}  = \zeta \left( r \right)\zeta \left( {m + 1,{{\left\{ 1 \right\}}_{p - 1}}} \right) - \zeta \left( {m + 1,{{\left\{ 1 \right\}}_{p - 1}},2,{{\left\{ 1 \right\}}_{r - 2}}} \right).\tag{2.36}\]
\end{thm}
\pf Similarly as in the proof of Theorem 2.6, we consider the multiple integral
\begin{align*}
&\int\limits_0^1 {\frac{1}{{{t_1}}}} d{t_1} \cdots \int\limits_0^{{t_{r - 1}}} {\frac{1}{{{t_r}}}} d{t_r}\int\limits_0^{{t_r}} {\frac{{\ln \left( {1 - {t_{r + 1}}} \right)}}{{{t_{r + 1}}}}} d{t_{r + 1}}\int\limits_0^{{t_{r + 1}}} {\frac{1}{{{t_{r + 2}}}}} d{t_{r + 2}} \cdots \int\limits_0^{{t_{p - 1}}} {\frac{1}{{{t_p}}}} d{t_p}\int\limits_0^{{t_p}} {\frac{{{{\ln }^m}\left( {1 - {t_{p + 1}}} \right)}}{{{t_{p + 1}}}}} d{t_{p + 1}}\\
& = \int\limits_{0 < {t_{p + 1}} <  \cdots  < {t_1} < 1} {\frac{{\ln \left( {1 - {t_{r + 1}}} \right){{\ln }^m}\left( {1 - {t_{p + 1}}} \right)}}{{{t_1}{t_2} \cdots {t_{p + 1}}}}} d{t_1}d{t_2} \cdots d{t_{p + 1}}.\tag{2.37}
\end{align*}
Applying the change of variables
${t_i} \mapsto 1 - {t_{p + 2 - i}}\quad (i = 1,2, \cdots ,p+ 1)$
to the above multiple integral, we deduce the identity
\begin{align*}
&\int\limits_{0 < {t_{p + 1}} <  \cdots  < {t_1} < 1} {\frac{{\ln \left( {1 - {t_{r + 1}}} \right){{\ln }^m}\left( {1 - {t_{p + 1}}} \right)}}{{{t_1}{t_2} \cdots {t_{p + 1}}}}} d{t_1}d{t_2} \cdots d{t_{p + 1}}\\
= &\int\limits_{0 < {t_{p + 1}} <  \cdots  < {t_1} < 1} {\frac{{\ln \left( {{t_{p - r + 1}}} \right){{\ln }^m}\left( {{t_1}} \right)}}{{\left( {1 - {t_{p + 1}}} \right)\left( {1 - {t_p}} \right) \cdots \left( {1 - {t_1}} \right)}}} d{t_1}d{t_2} \cdots d{t_{p + 1}}\\
 =& \int\limits_0^1 {\frac{{{{\ln }^m}\left( {{t_1}} \right)}}{{1 - {t_1}}}} d{t_1}\int\limits_0^{{t_1}} {\frac{1}{{1 - {t_2}}}} d{t_2} \cdots \int\limits_0^{{t_{p - r - 1}}} {\frac{1}{{1 - {t_{p - r}}}}} d{t_{p - r}}\int\limits_0^{{t_{p - r}}} {\frac{{\ln \left( {{t_{p - r + 1}}} \right)}}{{1 - {t_{p - r + 1}}}}} d{t_{p - r + 1}}\\
&\quad \times \int\limits_0^{{t_{p - r + 1}}} {\frac{1}{{1 - {t_{p - r + 2}}}}} d{t_{p - r + 2}} \cdots \int\limits_0^{{t_p}} {\frac{1}{{1 - {t_{p + 1}}}}} d{t_{p + 1}},\\
= &\sum\limits_{{n_1},{n_2}, \cdots {n_{p + 1}} = 1}^\infty  {\frac{{\left( {m + 1} \right)!{{\left( { - 1} \right)}^{m + 1}}}}{{{n_1} \cdots \left( {{n_1} +  \cdots  + {n_p}} \right){{\left( {{n_1} +  \cdots  + {n_p} + {n_{p + 1}}} \right)}^{m + 2}}}}} \\
& - \sum\limits_{{n_1},{n_2}, \cdots {n_{p + 1}} = 1}^\infty  {\frac{{m!{{\left( { - 1} \right)}^m}}}{{{n_1} \cdots \left( {{n_1} +  \cdots  + {n_{p - 1}}} \right){{\left( {{n_1} +  \cdots  + {n_p}} \right)}^2}{{\left( {{n_1} +  \cdots  + {n_p} + {n_{p + 1}}} \right)}^{m + 1}}}}} \\
& -  \cdots \\
& - \sum\limits_{{n_1},{n_2}, \cdots {n_{p + 1}} = 1}^\infty  {\frac{{{{\left( { - 1} \right)}^m}m!}}{{{n_1} \cdots \left( {{n_1} +  \cdots  + {n_r}} \right){{\left( {{n_1} +  \cdots  + {n_{r + 1}}} \right)}^2} \cdots {{\left( {{n_1} +  \cdots  + {n_p} + {n_{p + 1}}} \right)}^{m + 1}}}}} \\
 =& {\left( { - 1} \right)^{m + 1}}\left( {m + 1} \right)!\zeta \left( {m + 2,{{\left\{ 1 \right\}}_p}} \right) + {\left( { - 1} \right)^{m + 1}}m!\sum\limits_{i = 0}^{p - r - 1} {\zeta \left( {m + 1,{{\left\{ 1 \right\}}_i},2,{{\left\{ 1 \right\}}_{p - 1 - i}}} \right)} .\tag{2.38}
\end{align*}
On the other hand, by using (2.21), we have
\begin{align*}
\int\limits_0^1 {\frac{1}{{{t_1}}}} d{t_1} \cdots \int\limits_0^{{t_{r - 1}}} {\frac{1}{{{t_r}}}} d{t_r}\int\limits_0^{{t_r}} {\frac{{\ln \left( {1 - {t_{r + 1}}} \right)}}{{{t_{r + 1}^{1-n}}}}} d{t_{r + 1}} =& \frac{{{{\left( { - 1} \right)}^r}}}{{r!}}\int\limits_0^1 {t_1^{n - 1}{{\ln }^r}\left( {{t_1}} \right)\ln \left( {1 - {t_1}} \right)d{t_1}} \\
 =&  - \frac{{{H_n}}}{{{n^{r + 1}}}} - \sum\limits_{j = 1}^r {\frac{{{\zeta _n}\left( {j + 1} \right) - \zeta \left( {j + 1} \right)}}{{{n^{r - j + 1}}}}}.\tag{2.39}
\end{align*}
Substituting (2.39) into (2.37) and taking into account formula (2.26), we immediately arrive at
\begin{align*}
&\int\limits_{0 < {t_{p + 1}} <  \cdots  < {t_1} < 1} {\frac{{\ln \left( {1 - {t_{r + 1}}} \right){{\ln }^m}\left( {1 - {t_{p + 1}}} \right)}}{{{t_1}{t_2} \cdots {t_{p + 1}}}}} d{t_1}d{t_2} \cdots d{t_{p + 1}}\\
& = {\left( { - 1} \right)^{m + 1}}m!\sum\limits_{n = 1}^\infty  {\frac{{s\left( {n,m} \right)}}{{n!{n^{p - r}}}}\left\{ {\frac{{{H_n}}}{{{n^{r + 1}}}} + \sum\limits_{j = 1}^r {\frac{{{\zeta _n}\left( {j + 1} \right) - \zeta \left( {j + 1} \right)}}{{{n^{r - j + 1}}}}} } \right\}} .\tag{2.40}
\end{align*}
Combining (2.38) and (2.40), we get
\begin{align*}
&\sum\limits_{n = 1}^\infty  {\frac{{s\left( {n,m} \right)}}{{n!{n^{p - r}}}}\left\{ {\frac{{{H_n}}}{{{n^{r + 1}}}} + \sum\limits_{j = 1}^r {\frac{{{\zeta _n}\left( {j + 1} \right) - \zeta \left( {j + 1} \right)}}{{{n^{r - j + 1}}}}} } \right\}}  \\&= \left( {m + 1} \right)\zeta \left( {m + 2,{{\left\{ 1 \right\}}_p}} \right) + \sum\limits_{i = 0}^{p - r - 1} {\zeta \left( {m + 1,{{\left\{ 1 \right\}}_i},2,{{\left\{ 1 \right\}}_{p - 1 - i}}} \right)}.\tag{2.41}
 \end{align*}
Applying (2.29) and (2.30) into (2.41) yields
\[\sum\limits_{j = 1}^r {\sum\limits_{n = 1}^\infty  {\frac{{s\left( {n,m} \right){\zeta _n}\left( {j + 1} \right)}}{{n!{n^{p - j + 1}}}}} }  = \sum\limits_{j = 1}^r {\left\{ {\zeta \left( {j + 1} \right)\zeta \left( {m + 1,{{\left\{ 1 \right\}}_{p - j}}} \right) - \zeta \left( {m + 1,{{\left\{ 1 \right\}}_{p - j}},2,{{\left\{ 1 \right\}}_{j - 1}}} \right)} \right\}}.\tag{2.42}\]
Replacing $r$ by $r-1$ in (2.42), then substituting it into (2.42), we obtain
\[\sum\limits_{n = 1}^\infty  {\frac{{s\left( {n,m} \right){\zeta _n}\left( {r + 1} \right)}}{{n!{n^{p - r + 1}}}}}  = \zeta \left( {r + 1} \right)\zeta \left( {m + 1,{{\left\{ 1 \right\}}_{p - r}}} \right) - \zeta \left( {m + 1,{{\left\{ 1 \right\}}_{p - r}},2,{{\left\{ 1 \right\}}_{r - 1}}} \right).\]
The proof of Theorem 2.7 is finished.\hfill$\square$
\begin{re} In fact, by considering the following multiple integral
\[\int\limits_0^1 {\frac{{{{\ln }^r}\left( {{t_1}} \right)}}{{1 - {t_1}}}} d{t_1}\int\limits_0^{{t_1}} {\frac{1}{{{t_2}}}} d{t_2} \cdots \int\limits_0^{{t_{p - 1}}} {\frac{1}{{{t_p}}}} d{t_p}\int\limits_0^{{t_p}} {\frac{{{{\ln }^m}\left( {1 - {t_{p + 1}}} \right)}}{{{t_{p + 1}}}}} d{t_{p + 1}},\]
and using the elementary integral identity
\[\int\limits_0^1 {\frac{{{x^n}{{\ln }^m}x}}{{1 - x}}} dx = {\left( { - 1} \right)^m}m!\left( {\zeta \left( {m + 1} \right) - {\zeta _n}\left( {m + 1} \right)} \right),\;m \in \N.\]
then combining (2.25), we also obtain (2.36). When $p=1$, then the above integral can be rewritten as
\begin{align*}
 &\int\limits_0^1 {\frac{{{{\ln }^m}x}}{{1 - x}}} \int\limits_0^x {\frac{{{{\ln }^r}\left( {1 - t} \right)}}{t}} dtdx \\& = {\left( { - 1} \right)^{r + m}}r!m!\sum\limits_{n = r}^\infty  {\frac{{s\left( {n,r} \right)}}{{n!n}}\left( {\zeta \left( {m + 1} \right) - {\zeta _n}\left( {m + 1} \right)} \right)},\\
 &=\int\limits_0^1 {\frac{{{{\ln }^r}x}}{{1 - x}}} \int\limits_0^x {\frac{{{{\ln }^m}\left( {1 - t} \right)}}{t}} dtdx \\
  &= {\left( { - 1} \right)^{m + r}}m!r!\sum\limits_{n = r}^\infty  {\frac{{s\left( {n,m} \right)}}{{n!n}}\left( {\zeta \left( {r + 1} \right) - {\zeta _n}\left( {r + 1} \right)} \right)}.
 \end{align*}
 Applying (2.29), we deduce the following duality relation
\begin{align*}
 & \sum\limits_{n = 1}^\infty  {\frac{{s\left( {n,m} \right)}}{{n!n}}{\zeta _n}\left( {r + 1} \right)}  = \sum\limits_{n = 1}^\infty  {\frac{{s\left( {n,r} \right)}}{{n!n}}\zeta_n \left( {m + 1} \right)} ,
  \end{align*}
with $s\left( {n,k} \right) = 0,n < k,s\left( {n,0} \right) = s\left( {0,k} \right) = 0,s\left( {0,0} \right) = 1$. \hfill$\square$
\end{re}
In the same manner, we obtain the more general identity, see the following Theorem 2.9.
\begin{thm} For integers $k,m\in \N_0$ and $p\in \N$. Then
\[\sum\limits_{n = 1}^\infty  {\frac{{{Y_k}\left( n \right)s\left( {n,p} \right)}}{{n!{n^{m + 1}}}}}  = k!\sum\limits_{\scriptstyle 0 \le {i_1} +  \cdots +{i_m} \le k \hfill \atop
  \scriptstyle {i_j} \ge 0,j = 1,2 \cdots ,m \hfill} {\left( {\begin{array}{*{20}{c}}
   {p + k - \sum\limits_{j = 1}^m {{i_j}} }  \\
   p  \\
\end{array}} \right)\zeta \left( {p + k + 1 - \sum\limits_{j = 1}^m {{i_j}} ,{i_m} + 1, \cdots ,{i_1} + 1} \right)},\tag{2.43} \]
where $Y_k(n)$ is complete exponential Bell number defined by (2.8). When $m=0$, then
\begin{align*}
& \sum\limits_{n =1}^\infty  {\frac{{Y_k}\left( n \right){s\left( {n,p} \right)}}{{n!n}}}  = k!\left( {\begin{array}{*{20}{c}}
   {p + k}  \\
   k  \\
\end{array}} \right)\zeta \left( {p + k + 1} \right).
\end{align*}
\end{thm}
\pf In the same way as in proofs of Theorem 2.6 and 2.7, by using (2.7) and (2.26), we deduce that
\begin{align*}
&{\left( { - 1} \right)^{p + k}}p!\sum\limits_{n = 1}^\infty  {\frac{{{Y_k}\left( n \right)s\left( {n,p} \right)}}{{n!{n^{m + 1}}}}}  \\=& \int\limits_0^1 {\frac{{{{\ln }^k}\left( {1 - {t_1}} \right)}}{{{t_1}}}} d{t_1}\int\limits_0^{{t_1}} {\frac{1}{{{t_2}}}} d{t_2} \cdots \int\limits_0^{{t_{m - 1}}} {\frac{1}{{{t_m}}}} d{t_m}\int\limits_0^{{t_m}} {\frac{{{{\ln }^p}\left( {1 - {t_{m + 1}}} \right)}}{{{t_{m + 1}}}}} d{t_{m + 1}}\\
=& {\int _{0 < {t_{m + 1}} < {t_m} <  \cdots  < {t_1} < 1}}\frac{{{{\ln }^k}\left( {1 - {t_1}} \right){{\ln }^p}\left( {1 - {t_{m + 1}}} \right)}}{{{t_1}{t_2} \cdots {t_{m + 1}}}}d{t_1}d{t_2} \cdots d{t_{m + 1}}\\
=&\int\limits_{0 < {t_{m + 1}} <  \cdots  < {t_1} < 1} {\frac{{{{\ln }^p}\left( {{t_1}} \right){{\ln }^k}\left( {{t_{m + 1}}} \right)}}{{\left( {1 - {t_1}} \right) \cdots \left( {1 - {t_{m + 1}}} \right)}}} \\
 =&\int\limits_0^1 {\frac{{{{\ln }^p}\left( {{t_1}} \right)}}{{1 - {t_1}}}d{t_1}} \int\limits_0^{{t_1}} {\frac{1}{{1 - {t_2}}}d{t_2}}  \cdots \int\limits_0^{{t_{m - 1}}} {\frac{1}{{1 - {t_m}}}{d_m}} \int\limits_0^{{t_m}} {\frac{{{{\ln }^k}\left( {{t_{m + 1}}} \right)}}{{1 - {t_{m + 1}}}}d{t_{m + 1}}} \\
 =& \sum\limits_{{n_1}, \cdots ,{n_{m + 1}} = 1}^\infty  {\int\limits_0^1 {t_1^{{n_{m + 1}} - 1}{{\ln }^p}\left( {{t_1}} \right)d{t_1}} } \int\limits_0^{{t_1}} {t_2^{{n_m} - 1}d{t_2} \cdots \int\limits_0^{{t_{m - 1}}} {t_m^{{n_2} - 1}d{t_m}\int\limits_0^{{t_m}} {t_{m + 1}^{{n_1} - 1}{{\ln }^k}\left( {{t_{m + 1}}} \right)d{t_{m + 1}}} } }
.\tag{2.44}
\end{align*}
Then with the help of formula (2.33), we have
\begin{align*}
&p!\sum\limits_{n = 1}^\infty  {\frac{{{Y_k}\left( n \right)s\left( {n,p} \right)}}{{n!{n^{m + 1}}}}}\\
& = \sum\limits_{\scriptstyle 0 \le {i_1} +  \cdots +{i_m} \le k \hfill \atop
  \scriptstyle {i_j} \ge 0,j = 1,2 \cdots ,m \hfill} {{{\left( k \right)}_{{i_1}}}{{\left( {k - {i_1}} \right)}_{{i_2}}} \cdots {{\left( {k - \sum\limits_{j = 1}^{m - 1} {{i_j}} } \right)}_{{i_m}}}\left( {p + k - \sum\limits_{j = 1}^m {{i_j}} } \right)!} \\
  &\quad\quad\quad \times \zeta \left( {p + k + 1 - \sum\limits_{j = 1}^m {{i_j}} ,{i_m} + 1, \cdots ,{i_1} + 1} \right),\tag{2.45}
\end{align*}
where $(m)_l:=m(m-1)\cdots(m-l+1)$. By a direct calculation, we find that
\[{\left( k \right)_{{i_1}}}{\left( {k - {i_1}} \right)_{{i_2}}} \cdots {\left( {k - \sum\limits_{j = 1}^{m - 1} {{i_j}} } \right)_{{i_m}}}\left( {p + k - \sum\limits_{j = 1}^m {{i_j}} } \right)! = k!p!\left( {\begin{array}{*{20}{c}}
   {p + k - \sum\limits_{j = 1}^m {{i_j}} }  \\
   p  \\
\end{array}} \right).\tag{2.46}\]
Combining (2.45) and (2.46), we obtain the formula (2.43).\hfill$\square$\\
It is clear that Theorem 2.6 is an immediate corollary of Theorem 2.9.
Moreover, by a simple calculation, we also obtain the following result.
\begin{thm} For integers $k,m\in \N_0$ and $p\in \N$. Then
\begin{align*}
&p!\sum\limits_{n = 1}^\infty  {\frac{{{Y_k}\left( n \right)s\left( {n,p} \right)}}{{n!{n^{m + 1}}}}}  + {\left( { - 1} \right)^{m - 1}}k!\sum\limits_{n = 1}^\infty  {\frac{{{Y_p}\left( n \right)s\left( {n,k} \right)}}{{n!{n^{m + 1}}}}} \\
& = p!k!\sum\limits_{i = 1}^m {{{\left( { - 1} \right)}^{i - 1}}\zeta \left( {i + 1,{{\left\{ 1 \right\}}_{k - 1}}} \right)\zeta \left( {m + 2 - i,{{\left\{ 1 \right\}}_{p - 1}}} \right)}.\tag{2.47}
\end{align*}
\end{thm}
\pf First, by considering multiple integral, we can find that
\begin{align*}
\int\limits_0^1 {\frac{1}{{{t_1}}}d{t_1}}  \cdots \int\limits_0^{{t_{m - 1}}} {\frac{1}{{{t_m}}}} d{t_m}\int\limits_0^{{t_m}} {\frac{{{{\ln }^p}\left( {1 - {t_{m + 1}}} \right)}}{{{t_{m + 1}}}}d{t_{m + 1}}} & = {\left( { - 1} \right)^p}p!\sum\limits_{n = 1}^\infty  {\frac{{s\left( {n,p} \right)}}{{n!{n^{m + 1}}}}} \\
& = {\left( { - 1} \right)^p}p!\zeta \left( {m + 2,{{\left\{ 1 \right\}}_{p - 1}}} \right).\tag{2.48}
\end{align*}
Then applying formula (2.48) to the right hand side of (2.44) and integrating by parts, we show that
\begin{align*}
&{\left( { - 1} \right)^{p + k}}p!\sum\limits_{n = 1}^\infty  {\frac{{{Y_k}\left( n \right)s\left( {n,p} \right)}}{{n!{n^{m + 1}}}}} \\
 =& {\left( { - 1} \right)^{p + k}}p!k!\zeta \left( {2,{{\left\{ 1 \right\}}_{k - 1}}} \right)\zeta \left( {m + 1,{{\left\{ 1 \right\}}_{p - 1}}} \right)\\
& - \int\limits_0^1 {\frac{1}{{{t_1}}}d{t_1}\int\limits_0^{{t_1}} {\frac{{{{\ln }^k}\left( {1 - {t_2}} \right)}}{{{t_2}}}d{t_2}} } \int\limits_0^{{t_1}} {\frac{1}{{{t_2}}}} d{t_2} \cdots \int\limits_0^{{t_{m - 2}}} {\frac{1}{{{t_{m - 1}}}}} d{t_{m - 1}}\int\limits_0^{{t_{m - 1}}} {\frac{{{{\ln }^p}\left( {1 - {t_m}} \right)}}{{{t_m}}}d{t_m}} \\
=&  \cdots \\
=& {\left( { - 1} \right)^{p + k}}p!k!\sum\limits_{i = 1}^m {{{\left( { - 1} \right)}^{i - 1}}\zeta \left( {i + 1,{{\left\{ 1 \right\}}_{k - 1}}} \right)\zeta \left( {m + 2 - i,{{\left\{ 1 \right\}}_{p - 1}}} \right)} \\
 &+ {\left( { - 1} \right)^m}\int\limits_0^1 {\frac{{{{\ln }^p}\left( {1 - {t_1}} \right)}}{{{t_1}}}} d{t_1}\int\limits_0^{{t_1}} {\frac{1}{{{t_2}}}} d{t_2} \cdots \int\limits_0^{{t_{m - 1}}} {\frac{1}{{{t_m}}}} d{t_m}\int\limits_0^{{t_m}} {\frac{{{{\ln }^k}\left( {1 - {t_{m + 1}}} \right)}}{{{t_{m + 1}}}}} d{t_{m + 1}}\\
 = &{\left( { - 1} \right)^{p + k}}p!k!\sum\limits_{i = 1}^m {{{\left( { - 1} \right)}^{i - 1}}\zeta \left( {i + 1,{{\left\{ 1 \right\}}_{k - 1}}} \right)\zeta \left( {m + 2 - i,{{\left\{ 1 \right\}}_{p - 1}}} \right)} \\
 & + {\left( { - 1} \right)^{p + k + m}}k!\sum\limits_{n = 1}^\infty  {\frac{{{Y_p}\left( n \right)s\left( {n,k} \right)}}{{n!{n^{m + 1}}}}}.
\end{align*}
Thus, formula (2.47) holds.\hfill$\square$\\
Setting $m=1$ in the above equation we obtain
\[p!\sum\limits_{n = 1}^\infty  {\frac{{s\left( {n,p} \right){Y_k}\left( n \right)}}{{n!{n^2}}}}  + k!\sum\limits_{n = 1}^\infty  {\frac{{s\left( {n,k} \right){Y_p}\left( n \right)}}{{n!{n^2}}}}  = k!p!\zeta \left( {k + 1} \right)\zeta \left( {p + 1} \right).\tag{2.48}\]

\section{Relations between multiple zeta values and multiple zeta star values}
In this section, we will establish some explicit relationships between multiple zeta star values and multiple zeta
values. Furthermore, we give closed form for several classes of nonlinear Euler sums in terms of zeta values and linear sums.
\begin{lem} For positive integers $n$ and $p$, then
\begin{align*}
&\int\limits_0^1 {{x^{n - 1}}{\rm Li}{_p}\left( x \right)} dx = \sum\limits_{i = 1}^{p- 1} {\frac{{{{\left( { - 1} \right)}^{i - 1}}}}{{{n^i}}}} \zeta \left( {p + 1 - i} \right) + \frac{{{{\left( { - 1} \right)}^{p - 1}}}}{{{n^p}}}{H_n},\tag{3.1}\\
&\int\limits_0^1 {{x^{n - 1}}\ln x{\rm{L}}{{\rm{i}}_p}\left( x \right)} dx = \sum\limits_{i = 1}^{p - 1} {\sum\limits_{j = 1}^{p - i} {{{\left( { - 1} \right)}^{i + j - 1}}\frac{{\zeta \left( {p + 2 - i - j} \right)}}{{{n^{i + j}}}}} } \\
&\quad\quad\quad\quad\quad\quad\quad\quad\quad\quad+ {\left( { - 1} \right)^p}p\frac{{{H_n}}}{{{n^{p + 1}}}} + {\left( { - 1} \right)^p}\frac{{{\zeta _n}\left( 2 \right) - \zeta \left( 2 \right)}}{{{n^p}}}, \tag{3.2}
\end{align*}\
where ${\rm Li}{_p}\left( x \right)$ denotes the polylogarithm function defined by
\[{\rm Li}{_p}\left( x \right) := \sum\limits_{n = 1}^\infty  {\frac{{{x^n}}}{{{n^p}}}}, \Re(p)>1,\ \left| x \right| \le 1 ,\]
with ${\rm Li}{_1}\left( x \right)=-\ln(1-x),\ x\in [-1,1).$
\end{lem}
\pf Using integration by parts, we deduce that, for $x\in [-1,1)$,
\[\int\limits_0^x {{t^{n - 1}}{\rm Li}{_p}\left( t \right)dt}  = \sum\limits_{i = 1}^{p - 1} {{{\left( { - 1} \right)}^{i - 1}}\frac{{{x^n}}}{{{n^i}}}{\rm Li}{_{p + 1 - i}}\left( x \right)}  + \frac{{{{\left( { - 1} \right)}^p}}}{{{n^p}}}\ln \left( {1 - x} \right)\left( {{x^n} - 1} \right) - \frac{{{{\left( { - 1} \right)}^p}}}{{{n^p}}}\left( {\sum\limits_{k = 1}^n {\frac{{{x^k}}}{k}} } \right).\tag{3.3}\]
Letting $x$ approach 1, the result is (3.1). Multiplying (3.3) by $\frac{{{{\ln }^{m-1}}x}}{{ x}}$ and integrating over the interval $(0,1)$, we have the following recurrence relation
\begin{align*}
\int\limits_0^1 {{x^{n - 1}}{{\ln }^m}x{\rm{L}}{{\rm{i}}_p}\left( x \right)} dx =& m\sum\limits_{i = 1}^{p - 1} {\frac{{{{\left( { - 1} \right)}^i}}}{{{n^i}}}} \int\limits_0^1 {{x^{n - 1}}{{\ln }^{m - 1}}x{\rm{L}}{{\rm{i}}_{p + 1 - i}}\left( x \right)} dx + m!{\left( { - 1} \right)^{m + p - 1}}\frac{{{\zeta _n}\left( {m + 1} \right)}}{{{n^p}}}\\
&\quad + m!\frac{{{{\left( { - 1} \right)}^{m + p - 1}}}}{{{n^p}}}\left\{ {\frac{{{H_n}}}{{{n^m}}} - \sum\limits_{j = 1}^{m - 1} {\frac{{\zeta \left( {j + 1} \right) - {\zeta _n}\left( {j + 1} \right)}}{{{n^{m - j}}}}}  - \zeta \left( {m + 1} \right)} \right\}.\tag{3.4}
\end{align*}
Taking $m=1$ in (3.4), we obtain (3.2). This completes the proof of Lemma 3.1.\hfill$\square$
\begin{thm} For integers $p,m\in \N$, then we have the relation
\begin{align*}
{\zeta ^ \star }\left( {p + 1,{{\left\{ 1 \right\}}_m}} \right) =& \sum\limits_{i = 1}^{p - 1} {{{\left( { - 1} \right)}^{i - 1}}\zeta \left( {p + 1 - i} \right)\zeta \left( {m + 1,{{\left\{ 1 \right\}}_{i - 1}}} \right)} \\
&+ {\left( { - 1} \right)^{p + 1}}\sum\limits_{i = 1}^{p - 1} {\zeta \left( {m + 1,{{\left\{ 1 \right\}}_{i - 1}},2,{{\left\{ 1 \right\}}_{p - 1 - i}}} \right)} \\
& + {\left( { - 1} \right)^{p + 1}}\left( {m + 1} \right)\zeta \left( {m + 2,{{\left\{ 1 \right\}}_{p - 1}}} \right).\tag{3.5}
\end{align*}
\end{thm}
\pf First, we consider the integral
\begin{align*}
\int\limits_0^1 {\frac{{{\rm{L}}{{\rm{i}}_p}\left( x \right){{\ln }^m}\left( {1 - x} \right)}}{x}} dx =& \sum\limits_{n = 1}^\infty  {\frac{1}{{{n^p}}}\int\limits_0^1 {{x^{n - 1}}{{\ln }^m}\left( {1 - x} \right)} dx} \\
& = m!{(-1)^{m}}\sum\limits_{n = 1}^\infty  {\frac{{{\zeta ^ \star }\left( {{{\left\{ 1 \right\}}_m}} \right)}}{{{n^{p + 1}}}}} \\
& = m!{(-1)^{m}}{\zeta ^ \star }\left( {p + 1,{{\left\{ 1 \right\}}_m}} \right).\tag{3.6}
\end{align*}
On the other hand, by using (3.1), we conclude that
\begin{align*}
\int\limits_0^1 {\frac{{{\rm Li}{_p}\left( x \right){{\ln }^m}\left( {1 - x} \right)}}{x}} dx
 &={\left( { - 1} \right)^m}m!\sum\limits_{n = m}^\infty  {\frac{{s\left( {n,m} \right)}}{{n!}}\int\limits_0^1 {{x^{n - 1}}{\rm Li}{_p}\left( x \right)} dx}
  \nonumber \\
           &={\left( { - 1} \right)^m}m!\sum\limits_{i = 1}^{p - 1} {{{\left( { - 1} \right)}^{i - 1}}\zeta \left( {p + 1 - i} \right)\sum\limits_{n = m}^\infty  {\frac{{s\left( {n,m} \right)}}{{n!{n^i}}}} }
  \nonumber \\
           & \quad + {\left( { - 1} \right)^{m + p - 1}}m!\sum\limits_{n = 1}^\infty  {\frac{{{H_n}s\left( {n,m} \right)}}{{{n^p}n!}}} .\tag{3.7}
\end{align*}
Therefore, taking into account formulas (3.6), (3.7) and (2.29), we arrive at
\[{\zeta ^ \star }\left( {p + 1,{{\left\{ 1 \right\}}_m}} \right) = \sum\limits_{i = 0}^{p - 2} {{{\left( { - 1} \right)}^i}\zeta \left( {p - i} \right)\zeta \left( {m + 1,{{\left\{ 1 \right\}}_i}} \right)}  + {\left( { - 1} \right)^{p - 1}}\sum\limits_{n = 1}^\infty  {\frac{{{H_n}s\left( {n,m} \right)}}{{{n^p}n!}}}.\tag{3.8} \]
The formula (3.8) also immediately follows from (2.47). Then,
substituting (2.30) into (3.8) yields the desired result. \hfill$\square$
\begin{thm} For integers $p,m\in \N$, then we have the relation
\begin{align*}
&\sum\limits_{i = 1}^m {{\zeta ^ \star }\left( {p + 1,{{\left\{ 1 \right\}}_{i - 1}},2,{{\left\{ 1 \right\}}_{m - i}}} \right)}  + {\zeta ^ \star }\left( {p + 2,{{\left\{ 1 \right\}}_m}} \right) - \sum\limits_{i = 0}^{m - 1} {\zeta \left( {m - i + 1} \right){\zeta ^ \star }\left( {p + 1,{{\left\{ 1 \right\}}_i}} \right)} \\
& = \sum\limits_{i = 1}^{p - 1} {\sum\limits_{j = 1}^{p - i} {{{\left( { - 1} \right)}^{i + j}}\zeta \left( {p + 2 - i - j} \right)\zeta \left( {m + 1,{{\left\{ 1 \right\}}_{i + j - 1}}} \right)} }  + {\left( { - 1} \right)^{p + 1}}p\left( {m + 1} \right)\zeta \left( {m + 2,{{\left\{ 1 \right\}}_p}} \right)\\
&\quad + {\left( { - 1} \right)^{p + 1}}p\sum\limits_{i = 0}^{p - 1} {\zeta \left( {m + 1,{{\left\{ 1 \right\}}_i},2,{{\left\{ 1 \right\}}_{p - 1 - i}}} \right)}  - {\left( { - 1} \right)^{p + 1}}\zeta \left( {m + 1,{{\left\{ 1 \right\}}_{p - 1}},2} \right).\tag{3.9}
\end{align*}
\end{thm}
\pf Similarly as in the proof of Theorem 3.2, first, by using (2.6), we obtain the following identity
\begin{align*}
\int\limits_0^1 {\frac{{{\rm{L}}{{\rm{i}}_p}\left( x \right)\ln x{{\ln }^m}\left( {1 - x} \right)}}{x}} dx =& \sum\limits_{n = 1}^\infty  {\frac{1}{{{n^p}}}\int\limits_0^1 {{x^{n - 1}}\ln x{{\ln }^m}\left( {1 - x} \right)} dx} \\
= &m!{\left( { - 1} \right)^{m - 1}}\sum\limits_{i = 1}^m {{\zeta ^ \star }\left( {p + 1,{{\left\{ 1 \right\}}_{i - 1}},2,{{\left\{ 1 \right\}}_{m - i}}} \right)}\\
&+ m!{\left( { - 1} \right)^{m - 1}}{\zeta ^ \star }\left( {p + 2,{{\left\{ 1 \right\}}_m}} \right)\\
& - m!{\left( { - 1} \right)^{m - 1}}\sum\limits_{i = 0}^{m - 1} {\zeta \left( {m - i + 1} \right){\zeta ^ \star }\left( {p + 1,{{\left\{ 1 \right\}}_i}} \right)}.\tag{3.10}
\end{align*}
Then by using (2.26) and (3.2), we deduce the result
\begin{align*}
&\int\limits_0^1 {\frac{{{\rm{L}}{{\rm{i}}_p}\left( x \right)\ln x{{\ln }^m}\left( {1 - x} \right)}}{x}} dx\\ =&{\left( { - 1} \right)^m}m!\sum\limits_{n = 1}^\infty  {\frac{{s\left( {n,m} \right)}}{{n!}}} \int\limits_0^1 {{x^{n - 1}}{\rm{ln}}x{\rm{L}}{{\rm{i}}_p}\left( x \right)} dx\\
=&{\left( { - 1} \right)^m}m!\sum\limits_{i = 1}^{p - 1} {\sum\limits_{j = 1}^{p - i} {{{\left( { - 1} \right)}^{i + j-1}}\zeta \left( {p + 2 - i - j} \right)\zeta \left( {m + 1,{{\left\{ 1 \right\}}_{i + j - 1}}} \right)} } \\
 &+ {\left( { - 1} \right)^{m + p}}pm!\sum\limits_{n = 1}^\infty  {\frac{{{H_n}s\left( {n,m} \right)}}{{n!{n^{p + 1}}}}}  + {\left( { - 1} \right)^{m + p}}m!\sum\limits_{n = 1}^\infty  {\frac{{s\left( {n,m} \right)\left( {{\zeta _n}\left( 2 \right) - \zeta \left( 2 \right)} \right)}}{{n!{n^{p + 1}}}}} .\tag{3.11}
\end{align*}
Combining (2.30), (2.36), (3.10) and (3.11), we obtain (3.9). \hfill$\square$\\
Putting $p=1$ in (3.9), we can give the following Corollary.
\begin{cor} For positive integer $m$, we have
\begin{align*}
&\sum\limits_{i = 1}^m {{\zeta ^ \star }\left( {2,{{\left\{ 1 \right\}}_{i - 1}},2,{{\left\{ 1 \right\}}_{m - i}}} \right)}  + {\zeta ^ \star }\left( {3,{{\left\{ 1 \right\}}_m}} \right)\\ &= \sum\limits_{i = 0}^{m - 1} {\zeta \left( {m - i + 1} \right){\zeta ^ \star }\left( {2,{{\left\{ 1 \right\}}_i}} \right)}  + \left( {m + 1} \right)\zeta \left( {m + 2,1} \right).\tag{3.12}
\end{align*}
\end{cor}
Setting $p=1,2$ in (3.5) we obtain
\[\zeta^\star \left( {2,{{\left\{ 1 \right\}}_m}} \right) = \left( {m + 1} \right)\zeta \left( {m + 2} \right),\tag{3.13}\]
\[{\zeta ^ \star }\left( {3,{{\left\{ 1 \right\}}_m}} \right) = \zeta \left( 2 \right)\zeta \left( {m + 1} \right) - \left( {m + 1} \right)\zeta \left( {m + 2,1} \right) - \zeta \left( {m + 1,2} \right).\tag{3.14}\]
At the end of this section we give some closed form for several classes of nonlinear Euler sums. In \cite{BG1996}, J.M. Borwein, R. Girgensohn proved that all $\zeta \left( {q,p,r} \right)$ with $r+p+q$ is even or less than or equal to 10 or $r+p+q=12$ were reducible to zeta values and linear sums. From \cite{BG1996}, we have
\begin{align*}
&\zeta \left( {5,1} \right) = \frac{3}{4}\zeta \left( 6 \right) - \frac{1}{2}{\zeta ^2}\left( 3 \right),\\
&\zeta \left( {6,1} \right) = 3\zeta \left( 7 \right) - \zeta \left( 2 \right)\zeta \left( 5 \right) - \zeta \left( 3 \right)\zeta \left( 4 \right),\\
&\zeta \left( {6,1,1} \right) = \frac{{61}}{{24}}\zeta \left( 8 \right) - 3\zeta \left( 3 \right)\zeta \left( 5 \right) + \frac{1}{2}\zeta \left( 2 \right){\zeta ^2}\left( 3 \right),\\
&\zeta \left( {5,1,1,1} \right) = \frac{{499}}{{192}}\zeta \left( 8 \right) - 4\zeta \left( 3 \right)\zeta \left( 5 \right) + \zeta \left( 2 \right){\zeta ^2}\left( 3 \right),\\
&\zeta \left( {5,1,2} \right) =  - \frac{{73}}{{72}}\zeta \left( 8 \right) + \frac{9}{2}\zeta \left( 3 \right)\zeta \left( 5 \right) - \frac{3}{2}\zeta \left( 2 \right){\zeta ^2}\left( 3 \right) - {S_{2,6}},\\
&\zeta \left( {5,2,1} \right) =  - \frac{{541}}{{144}}\zeta \left( 8 \right) + \frac{7}{2}\zeta \left( 3 \right)\zeta \left( 5 \right) - \zeta \left( 2 \right){\zeta ^2}\left( 3 \right) + \frac{7}{2}{S_{2,6}},\\
&\zeta \left( {7,1,1} \right) = \frac{{28}}{3}\zeta \left( 9 \right) - 3\zeta \left( 2 \right)\zeta \left( 7 \right) - \frac{7}{4}\zeta \left( 3 \right)\zeta \left( 6 \right) - \frac{9}{4}\zeta \left( 4 \right)\zeta \left( 5 \right) + \frac{1}{6}{\zeta ^3}\left( 3 \right),\\
&\zeta \left( {6,1,2} \right) =  - \frac{{313}}{{36}}\zeta \left( 9 \right) + 7\zeta \left( 2 \right)\zeta \left( 7 \right) - \frac{5}{3}\zeta \left( 3 \right)\zeta \left( 6 \right) - \frac{1}{4}\zeta \left( 4 \right)\zeta \left( 5 \right) - \frac{1}{3}{\zeta ^3}\left( 3 \right),\\
&\zeta \left( {6,2,1} \right) =  - \frac{{2189}}{{72}}\zeta \left( 9 \right) + 11\zeta \left( 2 \right)\zeta \left( 7 \right) + \frac{9}{2}\zeta \left( 3 \right)\zeta \left( 6 \right) + \frac{{13}}{2}\zeta \left( 4 \right)\zeta \left( 5 \right) - \frac{1}{3}{\zeta ^3}\left( 3 \right).
\end{align*}
Putting $p=3,m=3$ in (3.5), we deduce that
\[\zeta^\star\left( {5,1,1,1} \right) =  - \frac{{385}}{{192}}\zeta \left( 8 \right) + 5\zeta \left( 3 \right)\zeta \left( 5 \right) - \zeta \left( 2 \right){\zeta ^2}\left( 3 \right) - \frac{3}{4}{S_{2,6}}.\]
Taking $p=3,m=4$ and $5$ in (3.5), we get
\begin{align*}
{\zeta ^ \star }\left( {4,{{\left\{ 1 \right\}}_4}} \right) =& \zeta \left( 3 \right)\zeta \left( 5 \right) - \zeta \left( 2 \right)\zeta \left( {5,1} \right) + 5\zeta \left( {6,1,1} \right) + \zeta \left( {5,2,1} \right) + \zeta \left( {5,1,2} \right)\\
=& \frac{{107}}{{16}}\zeta \left( 8 \right) - 6\zeta \left( 3 \right)\zeta \left( 5 \right) + \frac{1}{2}\zeta \left( 2 \right){\zeta ^2}\left( 3 \right) + \frac{3}{4}{S_{2,6}},\tag{3.15}\\
{\zeta ^ \star }\left( {4,{{\left\{ 1 \right\}}_5}} \right) =& \zeta \left( 3 \right)\zeta \left( 6 \right) - \zeta \left( 2 \right)\zeta \left( {6,1} \right) + 6\zeta \left( {7,1,1} \right) + \zeta \left( {6,2,1} \right) + \zeta \left( {6,1,2} \right)\\
=&\frac{{1217}}{{72}}\zeta \left( 9 \right) - 3\zeta \left( 2 \right)\zeta \left( 7 \right) - \frac{{59}}{{12}}\zeta \left( 3 \right)\zeta \left( 6 \right) - \frac{{19}}{4}\zeta \left( 4 \right)\zeta \left( 5 \right) + \frac{1}{3}{\zeta ^3}\left( 3 \right).
\end{align*}
where ${S_{2,6}} = {\zeta ^ \star }\left( {6,2} \right)$. On the other hand, we note that
\begin{align*}
{\zeta ^ \star }\left( {4,{{\left\{ 1 \right\}}_4}} \right) =& \sum\limits_{n = 1}^\infty  {\frac{{\zeta _n^ \star \left( {{{\left\{ 1 \right\}}_4}} \right)}}{{{n^4}}}}  = \frac{1}{{24}}\sum\limits_{n = 1}^\infty  {\frac{{{Y_4}\left( n \right)}}{{{n^4}}}} \\
=&\frac{1}{{24}}\sum\limits_{n = 1}^\infty  {\frac{{H_n^4 + 8{H_n}{\zeta _n}\left( 3 \right) + 6H_n^2{\zeta _n}\left( 2 \right) + 3\zeta _n^2\left( 2 \right) + 6{\zeta _n}\left( 4 \right)}}{{{n^4}}}}.\tag{3.16}
\end{align*}
In \cite{X2016}, we proved the results
\begin{align*}
&S_{13,4}=\sum\limits_{n = 1}^\infty  {\frac{{{H_n}{\zeta _n}\left( 3 \right)}}{{{n^4}}}}  =  - \frac{{511}}{{144}}\zeta \left( 8 \right) + 7\zeta \left( 3 \right)\zeta \left( 5 \right) + \zeta \left( 2 \right){\zeta ^2}\left( 3 \right) - \frac{{25}}{4}S_{2,6},\tag{3.17}\\
&S_{2^2,4}=\sum\limits_{n = 1}^\infty  {\frac{{\zeta _n^2\left( 2 \right)}}{{{n^4}}}}  = 11{S_{2,6}} + \frac{{457}}{{18}}\zeta \left( 8 \right) + 6\zeta \left( 2 \right){\zeta ^2}\left( 3 \right) - 40\zeta \left( 3 \right)\zeta \left( 5 \right),\tag{3.18}
\end{align*}
and
\[\sum\limits_{n = 1}^\infty  {\frac{{6H_n^2{\zeta _n}\left( 2 \right) - H_n^4}}{{{n^4}}}}  =  - 17S_{2,6} - 10\zeta \left( 2 \right){\zeta ^2}\left( 3 \right) + 104\zeta \left( 3 \right)\zeta \left( 5 \right) - \frac{{5911}}{{72}}\zeta \left( 8 \right),\tag{3.19}\]
Substituting (3.17) and (3.18) into (3.16) respectively, we have
\[\sum\limits_{n = 1}^\infty  {\frac{{H_n^4 + 6H_n^2{\zeta _n}\left( 2 \right)}}{{{n^4}}}}  = \frac{{956}}{9}\zeta \left( 8 \right) - 80\zeta \left( 3 \right)\zeta \left( 5 \right) - 14\zeta \left( 2 \right){\zeta ^2}\left( 3 \right) + 35{S_{2,6}}.\tag{3.20}\]
Hence, the relations (3.19) and (3.20) yield the following identities
\begin{align*}
&S_{1^4,4}=\sum\limits_{n = 1}^\infty  {\frac{{H_n^4}}{{{n^4}}}}  = \frac{{13559}}{{144}}\zeta \left( 8 \right) - 92\zeta \left( 3 \right)\zeta \left( 5 \right) - 2\zeta \left( 2 \right){\zeta ^2}\left( 3 \right) + 26{S_{2,6}},\tag{3.21}\\
&S_{1^22,4}=\sum\limits_{n = 1}^\infty  {\frac{{H_n^2{\zeta _n}\left( 2 \right)}}{{{n^4}}}}  = \frac{{193}}{{96}}\zeta \left( 8 \right) + 2\zeta \left( 3 \right)\zeta \left( 5 \right) - 2\zeta \left( 2 \right){\zeta ^2}\left( 3 \right) + \frac{3}{2}{S_{2,6}}.\tag{3.22}
\end{align*}
In fact, using the method of this paper, it is possible to evaluate other
Euler sums involving harmonic numbers. For example, we have used our method to obtain the following explicit evaluations. (Some of them were given previously by David H. Bailey, Jonathan M. Borwein, Roland Girgensohn \cite{BBG1994} and Philippe Flajolet, Bruno Salvy\cite{FS1998})
\begin{exa} Some examples on nonlinear Euler sums follows:
\begin{align*}
&S_{1^22,2}=\sum\limits_{n = 1}^\infty  {\frac{{H_n^2{\zeta _n}(2)}}{{{n^2}}}}  = \frac{{41}}{{12}}\zeta \left( 6 \right) + 2{\zeta ^2}\left( 3 \right),\\
&S_{1^22,3}=\sum\limits_{n = 1}^\infty  {\frac{{H_n^2{\zeta _n}\left( 2 \right)}}{{{n^3}}}}  =  - 7\zeta \left( 7 \right) + \frac{{19}}{2}\zeta \left( 3 \right)\zeta \left( 4 \right) - 2\zeta \left( 2 \right)\zeta \left( 5 \right),\\
&S_{1^32,3}=\sum\limits_{n =
1}^\infty {\frac{{H_n^3{\zeta _n}(2)}}{{{n^2}}}}  =
\frac{{83}}{{16}}\zeta (7) + \frac{{27}}{2}\zeta (3)\zeta (4)- \frac{5}{2}\zeta (2)\zeta (5) ,\\
&S_{1^22,3}=\sum\limits_{n = 1}^\infty {\frac{{H_n^2{\zeta _n}(2)}}{{{n^3}}}}
=  - 7\zeta (7) + \frac{{19}}{2}\zeta (3)\zeta (4) -
2\zeta (2)\zeta (5),\\
&S_{1^23,2}=\sum\limits_{n = 1}^\infty {\frac{{H_n^2{\zeta _n}(3)}}{{{n^2}}}}
= \frac{{329}}{{16}}\zeta (7)- 6\zeta (3)\zeta (4) -
\frac{9}{2}\zeta
(2)\zeta (5),\\
&S_{12^2,2}=\sum\limits_{n = 1}^\infty  {\frac{{{H_n}\zeta
_n^2(2)}}{{{n^2}}}}  =  - \frac{{217}}{{16}}\zeta (7) + 5\zeta
(3)\zeta (4) + \frac{{13}}{2}\zeta (2)\zeta (5),\\
&S_{1^24,2}=\sum\limits_{n = 1}^\infty  {\frac{{H_n^2{\zeta _n}\left( 4 \right)}}{{{n^2}}}}  = \frac{{1289}}{{96}}\zeta \left( 8 \right) - 11\zeta \left( 3 \right)\zeta \left( 5 \right) + 5{S_{2,6}},\\
&S_{1^23,3}=\sum\limits_{n = 1}^\infty  {\frac{{H_n^2{\zeta _n}\left( 3 \right)}}{{{n^3}}}}  =  - \frac{{443}}{{288}}\zeta \left( 8 \right) + \frac{9}{2}\zeta \left( 3 \right)\zeta \left( 5 \right) + \frac{3}{2}\zeta \left( 2 \right){\zeta ^2}\left( 3 \right) - \frac{{23}}{4}{S_{2,6}},\\
&{S_{{1^2}{2^2},2}} = \sum\limits_{n = 1}^\infty  {\frac{{H_n^2\zeta _n^2\left( 2 \right)}}{{{n^2}}}}  = \frac{{55}}{8}\zeta \left( 8 \right) - 7\zeta \left( 3 \right)\zeta \left( 5 \right) + 2\zeta \left( 2 \right){\zeta ^2}\left( 3 \right) + 6{S_{2,6}},\\
&{S_{{3^2},2}} = \sum\limits_{n = 1}^\infty  {\frac{{\zeta _n^2\left( 3 \right)}}{{{n^2}}}}  = \frac{{677}}{{24}}\zeta \left( 8 \right) - 35\zeta \left( 3 \right)\zeta \left( 5 \right) + 4\zeta \left( 2 \right){\zeta ^2}\left( 3 \right) + \frac{{15}}{2}{S_{2,6}},\\
&{S_{23,3}} = \sum\limits_{n = 1}^\infty  {\frac{{{\zeta _n}\left( 2 \right){\zeta _n}\left( 3 \right)}}{{{n^3}}}}  =  - \frac{{827}}{{48}}\zeta \left( 8 \right)+\frac{{45}}{2}\zeta \left( 3 \right)\zeta \left( 5 \right) - \frac{3}{2}\zeta \left( 2 \right){\zeta ^2}\left( 3 \right) - \frac{{23}}{4}{S_{2,6}},\\
&{S_{24,2}} = \sum\limits_{n = 1}^\infty  {\frac{{\zeta _n^{}\left( 2 \right)\zeta _n^{}\left( 4 \right)}}{{{n^2}}}}  =  - \frac{{403}}{{36}}\zeta \left( 8 \right) + 20\zeta \left( 3 \right)\zeta \left( 5 \right) - 3\zeta \left( 2 \right){\zeta ^2}\left( 3 \right) - \frac{9}{2}{S_{2,6}},\\
&S_{1^5,3}=\sum\limits_{n = 1}^\infty  {\frac{{H_n^5}}{{{n^3}}}}  = \frac{{60499}}{{288}}\zeta \left( 8 \right) - \frac{{393}}{2}\zeta \left( 3 \right)\zeta \left( 5 \right) - \frac{{{\rm{15}}}}{{\rm{2}}}\zeta \left( 2 \right){\zeta ^2}\left( 3 \right) + \frac{{{\rm{235}}}}{{\rm{4}}}{S_{2,6}},\\
&{S_{{1^6},2}} = \sum\limits_{n = 1}^\infty  {\frac{{H_n^6}}{{{n^2}}}}  = \frac{{27903}}{6}\zeta \left( 8 \right) + 31\zeta \left( 3 \right)\zeta \left( 5 \right) + 16\zeta \left( 2 \right){\zeta ^2}\left( 3 \right) + 57{S_{2,6}},\\
&S_{1^32,3}=\sum\limits_{n = 1}^\infty  {\frac{{H_n^3{\zeta _n}\left( 2 \right)}}{{{n^3}}}}  =  - \frac{{2159}}{{48}}\zeta \left( 8 \right) + \frac{{93}}{2}\zeta \left( 3 \right)\zeta \left( 5 \right) + \frac{3}{2}\zeta \left( 2 \right){\zeta ^2}\left( 3 \right) - \frac{{53}}{4}{S_{2,6}},\\
&S_{12^2,3}=\sum\limits_{n = 1}^\infty  {\frac{{{H_n}\zeta _n^2\left( 2 \right)}}{{{n^3}}}}  =  - \frac{{6313}}{{288}}\zeta \left( 8 \right) + \frac{{43}}{2}\zeta \left( 3 \right)\zeta \left( 5 \right) + \frac{1}{2}\zeta \left( 2 \right){\zeta ^2}\left( 3 \right) - \frac{{17}}{4}{S_{2,6}},\\
&{S_{{1^4}2,2}} = \sum\limits_{n = 1}^\infty  {\frac{{H_n^4{\zeta _n}\left( 2 \right)}}{{{n^2}}}}  =  - \frac{{6631}}{{288}}\zeta \left( 8 \right) + 90\zeta \left( 3 \right)\zeta \left( 5 \right) + 3\zeta \left( 2 \right){\zeta ^2}\left( 3 \right) - \frac{{47}}{2}{S_{2,6}},\\
&{S_{123,2}} = \sum\limits_{n = 1}^\infty  {\frac{{{H_n}{\zeta _n}\left( 2 \right){\zeta _n}\left( 3 \right)}}{{{n^2}}}}  =  - \frac{{181}}{{288}}\zeta \left( 8 \right) + \frac{{15}}{2}\zeta \left( 3 \right)\zeta \left( 5 \right) - \frac{3}{2}\zeta \left( 2 \right){\zeta ^2}\left( 3 \right) - \frac{7}{4}{S_{2,6}},\\
&{S_{{1^3}3,2}} = \sum\limits_{n = 1}^\infty  {\frac{{H_n^3{\zeta _n}\left( 3 \right)}}{{{n^2}}}}  = \frac{{809}}{{48}}\zeta \left( 8 \right) + \frac{{23}}{2}\zeta \left( 3 \right)\zeta \left( 5 \right) - \frac{7}{2}\zeta \left( 2 \right){\zeta ^2}\left( 3 \right) - \frac{{33}}{4}{S_{2,6}},\\
&{S_{{1^5},4}} = \sum\limits_{n = 1}^\infty  {\frac{{H_n^5}}{{{n^4}}}}  = \frac{{4721}}{{36}}\zeta \left( 9 \right) + \frac{{265}}{8}\zeta \left( 2 \right)\zeta \left( 7 \right) - \frac{{4895}}{{24}}\zeta \left( 3 \right)\zeta \left( 6 \right) + 66\zeta \left( 4 \right)\zeta \left( 5 \right) - 5{\zeta ^3}\left( 3 \right),\\
&{S_{{1^2}3,4}} = \sum\limits_{n = 1}^\infty  {\frac{{H_n^2{\zeta _n}\left( 3 \right)}}{{{n^4}}}}  = \frac{{3895}}{{72}}\zeta \left( 9 \right) - \frac{5}{8}\zeta \left( 2 \right)\zeta \left( 7 \right) - \frac{{227}}{{24}}\zeta \left( 3 \right)\zeta \left( 6 \right) - \frac{{75}}{2}\zeta \left( 4 \right)\zeta \left( 5 \right) + {\zeta ^3}\left( 3 \right),\\
&{S_{{1^3}2,4}} = \sum\limits_{n = 1}^\infty  {\frac{{H_n^3{\zeta _n}\left( 2 \right)}}{{{n^4}}}}  =  - \frac{{449}}{{36}}\zeta \left( 9 \right) - 7\zeta \left( 2 \right)\zeta \left( 7 \right) + \frac{{11}}{8}\zeta \left( 3 \right)\zeta \left( 6 \right) + 27\zeta \left( 4 \right)\zeta \left( 5 \right) - \frac{{11}}{3}{\zeta ^3}\left( 3 \right),\\
&{S_{{{12}^2},4}} = \sum\limits_{n = 1}^\infty  {\frac{{{H_n}\zeta _n^2\left( 2 \right)}}{{{n^4}}}}  =  - \frac{{775}}{{36}}\zeta \left( 9 \right) + \frac{{85}}{8}\zeta \left( 2 \right)\zeta \left( 7 \right) - \frac{{221}}{{24}}\zeta \left( 3 \right)\zeta \left( 6 \right) + 10\zeta \left( 4 \right)\zeta \left( 5 \right) + 3{\zeta ^3}\left( 3 \right),\\
&S_{1^22,5}=\sum\limits_{n = 1}^\infty  {\frac{{H_n^2{\zeta _n}\left( 2 \right)}}{{{n^5}}}}  =  - \frac{{1481}}{{72}}\zeta \left( 9 \right) - 3{\zeta ^3}\left( 3 \right) - 5\zeta \left( 2 \right)\zeta \left( 7 \right) + \frac{{295}}{{24}}\zeta \left( 3 \right)\zeta \left( 6 \right) + 18\zeta \left( 4 \right)\zeta \left( 5 \right),\\
&{S_{{{13}^2},3}} = \sum\limits_{n = 1}^\infty  {\frac{{{H_n}\zeta _n^2\left( 3 \right)}}{{{n^3}}}}  = \frac{{883}}{{20}}\zeta \left( {10} \right) - 26{\zeta ^2}\left( 5 \right) - \frac{{31}}{4}\zeta \left( 3 \right)\zeta \left( 7 \right) - 8\zeta \left( 2 \right)\zeta \left( 3 \right)\zeta \left( 5 \right)\\
&\quad\quad\quad\quad\quad\quad\quad\quad\quad\quad\quad+ \frac{3}{4}{\zeta ^2}\left( 3 \right)\zeta \left( 4 \right) + 9\zeta \left( 2 \right){S_{2,6}} - \frac{{21}}{4}{S_{2,8}}.
\end{align*}
\end{exa}
Therefore, from Example 3.5 and references \cite{FS1998,X2016,X2017}, we obtain the following description of nonlinear Euler sums of weight $\leq 10$.
\begin{thm} All Euler sums of weight $\leq 8$ are reducible to $\mathbb{Q}$-linear combinations of single zeta monomials with the addition of $\{S_{2,6}\}$ for weight 8. For weight 9, all Euler sums of the form ${S_{{s_1} \cdots {s_k},q}}$ with $q\in \{4,5,6,7\}$ are expressible polynomially in terms of zeta values. For weight $p_1+p_2+q=10$, all quadratic sums $S_{p_1p_2,q}$ are reducible to $S_{2,6}$ and $S_{2,8}$.
\end{thm}
Moreover, we use Mathematica tool to check numerically each of the specific identities listed. The numerical values of
nonlinear Euler sums of weights 8 and 9, to 30 decimal digits, are:
\newcommand{\tabincell}[2]{\begin{tabular}{@{}#1@{}}#2\end{tabular}}
\begin{table}[htbp]
 \centering
\caption{\label{tab:test}Numerical approximation}
 \begin{tabular}{|c|c|c|}
  \hline
 \tabincell{c}{Euler sum} &\tabincell{c}{ Numerical values of closed form\\ (30 decimal digits)}
  & \tabincell{c}{Numerical approximation of\\ Euler sum (30 decimal digits)}\\
  \hline
  $S_{1^4,4}$ &1.68625748775730579166360832694&1.68625748775730579166360833402 \\
 \hline
  $S_{1^22,4}$ &1.29068714089613697618140723840&1.29068714089613697618140723441 \\
\hline
$S_{1^32,3}$ &2.82229596096025461026149662829&2.82229596096025461026149661993 \\
\hline
$S_{1^23,3}$&1.75388174691782356380634371202& 1.75388174691782356380634370478\\
\hline
$S_{12^2,3}$& 1.63443098048025390280783629910 &1.63443098048025390280783629225\\
\hline
$S_{1^24,2}$&4.88040799023015427295866390307& 4.88040799023015427295866390372\\
\hline
$S_{2^2,4}$& 1.13642391274089928376327915373 &1.13642391274089928376327915559\\
\hline
$S_{1^5,3}$&8.20602621468401623725850548850 &8.20602621468401623725850551862\\
\hline
$S_{123,2}$&3.36374308381687640081618084070 &3.36374308381687640081618083742\\
\hline
$S_{1^42,2}$&72.1778863641121208246730431963 &72.1778863641121208246730431899\\
\hline
$S_{1^33,2}$&14.5074537674864815323431145949 &14.5074537674864815323431145933\\
\hline
$S_{1^6,2}$&1302.28271941001924714647587730 &1302.28271941001924714647587732\\
\hline
$S_{1^5,4}$&2.31083536190405961638953653685 &2.31083536190405961638953653376\\
\hline
$S_{1^23,4}$&1.25355563158689137948838467515 &1.25355563158689137948838467072\\
\hline
$S_{12^2,4}$&1.22503753401105341474879224535 &1.22503753401105341474879224098\\
\hline
$S_{1^32,4}$&1.50676526085085659032600904678 &1.50676526085085659032600904154\\
\hline
 \end{tabular}
\end{table}
\\
In \cite{KO2010}, Masanobu Kaneko and Yasuo Ohno proved that
\begin{align*}
&{\left( { - 1} \right)^k}{\zeta ^ \star }\left( {k + 1,{{\left\{ 1 \right\}}_n}} \right) - {\left( { - 1} \right)^n}{\zeta ^ \star }\left( {n + 1,{{\left\{ 1 \right\}}_k}} \right)\\
& = k\zeta \left( {k + 2,{{\left\{ 1 \right\}}_{n - 1}}} \right) - n\zeta \left( {n + 2,{{\left\{ 1 \right\}}_{k - 1}}} \right)\\
&\quad + {\left( { - 1} \right)^k}\sum\limits_{j = 0}^{k - 2} {{{\left( { - 1} \right)}^j}\zeta \left( {k - j} \right)\zeta \left( {n + 1,{{\left\{ 1 \right\}}_j}} \right)} \\
&\quad - {\left( { - 1} \right)^n}\sum\limits_{j = 0}^{n - 2} {{{\left( { - 1} \right)}^j}\zeta \left( {n - j} \right)\zeta \left( {k + 1,{{\left\{ 1 \right\}}_j}} \right)} .\tag{3.23}
\end{align*}
Combining (2.30), (3.5), (3.8) and (3.23), we immediately deduce that
\begin{align*}
&\sum\limits_{n = 1}^\infty  {\frac{{{H_n}s\left( {n,p} \right)}}{{n!{n^m}}}}  - \sum\limits_{n = 1}^\infty  {\frac{{{H_n}s\left( {n,m} \right)}}{{n!{n^p}}}}  = p\zeta \left( {p + 2,{{\left\{ 1 \right\}}_{m - 1}}} \right) - m\zeta \left( {m + 2,{{\left\{ 1 \right\}}_{p - 1}}} \right),\tag{3.24}\\
&\sum\limits_{i = 1}^{p - 1} {\zeta \left( {m + 1,{{\left\{ 1 \right\}}_{i - 1}},2,{{\left\{ 1 \right\}}_{p - 1 - i}}} \right)}  - \sum\limits_{i = 1}^{m - 1} {\zeta \left( {p + 1,{{\left\{ 1 \right\}}_{i - 1}},2,{{\left\{ 1 \right\}}_{m - 1 - i}}} \right)} \\
& = \zeta \left( {p + 2,{{\left\{ 1 \right\}}_{m - 1}}} \right) - \zeta \left( {m + 2,{{\left\{ 1 \right\}}_{p - 1}}} \right).\tag{3.25}
\end{align*}
\section{Some identities for $H\left( {a,b;m,p} \right)$ and ${H^ \star }\left( {a,b;m,p} \right) $}
In 2012, Zagier \cite{DZ2012} proved that the multiple zeta star values ${\zeta ^ \star \left( {{{\left\{ 2 \right\}}_a},3,{{\left\{ 2 \right\}}_b}} \right)} $ and multiple zeta values ${\zeta \left( {{{\left\{ 2 \right\}}_a},3,{{\left\{ 2 \right\}}_b}} \right)} $ are reducible to polynomials in zeta values, $a,b\in \N_0$, and gave explicit formulae.
In this section, we will prove that the sums $H\left( {a,b;m,p} \right)$ and ${H^ \star }\left( {a,b;m,p} \right) $ can be expressed in terms of the Riemann zeta values, where $a,b\in \N_0$, the sums $H\left( {a,b;m,p} \right)$ and ${H^ \star }\left( {a,b;m,p} \right) $ are defined by
\[H\left( {a,b;m,p} \right) := \sum\limits_{a + b = m - 1} {\zeta \left( {{{\left\{ p \right\}}_a},p + 1,{{\left\{ p \right\}}_b}} \right)}, m,p\in \N, \]
\[{H^ \star }\left( {a,b;m,p} \right) := \sum\limits_{a + b = m - 1} {{\zeta ^ \star }\left( {{{\left\{ p \right\}}_a},p + 1,{{\left\{ p \right\}}_b}} \right)}, m,p\in \N. \]
Moreover, we investigate some cases in which the sums $H\left( {a,b;m,p} \right)$ and ${H^ \star }\left( {a,b;m,p} \right) $ are described by the Riemann zeta values. As a result, we obtain various relations between the multiple zeta (star) values and the Riemann zeta values, which contain most of the known results and some new ones.
Let $m$ be a positive integer, we define parametric Hurwitz zeta function and parametric Hurwitz zeta star function  by
\[\zeta \left( {{s_1},{s_2}, \cdots ,{s_m};a + 1} \right): = \sum\limits_{1 \le {k_1} <  \cdots  < {k_m}} {\frac{1}{{{{\left( {{k_1} + a} \right)}^{{s_1}}}{{\left( {{k_1} + a} \right)}^{{s_2}}} \cdots {{\left( {{k_m} + a} \right)}^{{s_m}}}}}} ,\]
\[{\zeta ^ \star }\left( {{s_1},{s_2}, \cdots ,{s_m};a + 1} \right): = \sum\limits_{1 \le {k_1} \le  \cdots  \le {k_m}} {\frac{1}{{{{\left( {{k_1} + a} \right)}^{{s_1}}}{{\left( {{k_1} + a} \right)}^{{s_2}}} \cdots {{\left( {{k_m} + a} \right)}^{{s_m}}}}}} ,\]
where $\Re(s_1)>1,\ s_i\geq 1,\ a\neq -1,-2,\cdots.$
When $m=1$, then the parametric Hurwitz zeta function ( or parametric Hurwitz zeta star function ) reduces to the classical Hurwitz zeta function, which is defined by
\[\zeta \left( {s,a + 1} \right) := \sum\limits_{n = 1}^\infty  {\frac{1}{{{{\left( {n + a} \right)}^s}}}}.\]
Next, we state and prove our main result on sums $H\left( {a,b;m,p} \right)$ and ${H^ \star }\left( {a,b;m,p} \right) $.
\begin{thm}
Define two sequences ${A_m(n)}$ and ${B_m(n)}$ by
$$ A_m(n) = (m-1)!\underset{i = 0}{\overset{m - 1}{\sum}}\dfrac{A_i(n)}{i!}\underset{k = 1}{\overset{n}{\sum}}x_k^{m - i},A_0(n) = 1,\left( {{x_k} \in \mathbb{C},k = 1,2, \cdots ,n} \right),$$
$$ B_m(n) = \underset{k_1 = 1}{\overset{n}{\sum}}x_{k_1}\underset{k_2 = 1}{\overset{k_1}{\sum}}x_{k_2}\cdots\underset{k_m = 1}{\overset{k_{m-1}}{\sum}}x_{k_m},B_0(n) = 1,\left( {{x_k} \in \mathbb{C},k = 1,2, \cdots ,n} \right).$$
Then \[A_m(n) = m!B_m(n).\tag{4.1}\]
\end{thm}
\pf
By induction, we can tell that $A_m(n)$ and $B_m(n)$ are polynomials of degree $m$ with $n$ variables $x_1, x_2, \cdots x_n$, moreover, $B_m(n)$ has coefficient $1$ at $x_1^{s_1}x_2^{s_2}\cdots x_n^{s_n}$, where $\underset{k = 1}{\overset{n}{\sum}}s_k = m$, and $s_k \geq 0$. Suppose $c_{s_1, \cdots, s_n}$ is the coefficient of $A_m(n)$ at $x_1^{s_1}x_2^{s_2}\cdots x_n^{s_n}$, then in order to get the conclusion, all we need to do is to prove that $c_{s_1, \cdots, s_n} = m!$ by induction.
For convenience, we can suppose $c_{t_1, \cdots, t_n} = 0$ if there is a $t_i < 0$, otherwise $c_{t_1, \cdots, t_n}$ is the coefficient of $A_{t_1 +\cdots + t_n}(n)$ at $x_1^{t_1}x_2^{t_2}\cdots x_n^{t_n}$.
For $m = 0$, the conclusion holds.
If the conclusion holds for any $i \leq m - 1$. For $m$, by the recursion formula of $A_m(n)$, $c_{s_1, \cdots, s_n} = (m-1)!\underset{i = 0}{\overset{m - 1}{\sum}}\underset{k = 1}{\overset{n}{\sum}}\dfrac{c_{s_1,\cdots,s_k - (m - i), \cdots, s_n}}{i!}$, since for $0 \leq i \leq m-1$,  $\dfrac{A_i(n)}{i!}\underset{k = 1}{\overset{n}{\sum}}x_k^{m - i}$ at $x_1^{s_1}x_2^{s_2}\cdots x_n^{s_n} = x_1^{s_1}\cdots x_k^{s_k - (m - i)}\cdots x_n^{s_n} \cdot x_k^{m - i}$, where $1 \leq k \leq n$, provides coefficient $\dfrac{c_{s_1,\cdots,s_k - (m - i), \cdots, s_n}}{i!}$. So
\begin{align*}
{c_{{s_1}, \cdots ,{s_n}}} &= \left( {m - 1} \right)!\sum\limits_{k = 1}^n {\sum\limits_{i = 0}^{m - 1} {\frac{{{c_{{s_1}, \cdots ,{s_k} - (m - i), \cdots ,{s_n}}}}}{{i!}}} } \\
& = \left( {m - 1} \right)!\sum\limits_{k = 1}^n {\sum\limits_{{i_k} = m - {s_k}}^{m - 1} {\frac{{{c_{{s_1}, \cdots ,{s_k} - (m - i), \cdots ,{s_n}}}}}{{i!}}} } \\
& = \left( {m - 1} \right)!\sum\limits_{k = 1}^n {\sum\limits_{{i_k} = m - {s_k}}^{m - 1} 1 } \\
& = \left( {m - 1} \right)!\sum\limits_{k = 1}^n {\left( {\left( {m - 1} \right) - \left( {m - {s_k}} \right) + 1} \right)} \\
&= \left( {m - 1} \right)!\sum\limits_{k = 1}^n {{s_k}} \\
&= m!,
\end{align*}
since $c_{s_1,\cdots,s_k - (m - i_k), \cdots, s_n} = i_k!$ by the inductive assumption. So the conclusion holds by induction. So far we have completed the proof of Theorem 4.1. \hfill$\square$
\\
Let ${X_n}\left( m \right) := \sum\limits_{i = 1}^n {x_i^m} $ in Theorem 4.1, we deduce the following identities
\begin{align*}
 &{B_0}\left( n \right) = 1,{B_1}\left( n \right) = {X_n}\left( 1 \right),{B_2}\left( n \right) = \frac{{X_n^2\left( 1 \right) + {X_n}\left( 2 \right)}}{2}, \\
 &{B_3}\left( n \right) = \frac{{X_n^3\left( 1 \right) + 3{X_n}\left( 1 \right){X_n}\left( 2 \right) + 2{X_n}\left( 3 \right)}}{{3!}}, \\
 &{B_4}\left( n \right) = \frac{{X_n^4\left( 1 \right) + 8{X_n}\left( 1 \right){X_n}\left( 3 \right) + 3X_n^2\left( 2 \right) + 6X_n^2\left( 1 \right){X_n}\left( 2 \right) + 6{X_n}\left( 4 \right)}}{{4!}}, \\
 &{B_5}\left( n \right) = \frac{1}{{5!}}\left\{ \begin{array}{l}
 X_n^5\left( 1 \right) + 10X_n^3\left( 1 \right){X_n}\left( 2 \right) + 20X_n^2\left( 1 \right){X_n}\left( 3 \right) + 15{X_n}\left( 1 \right)X_n^2\left( 2 \right) \\
  + 30{X_n}\left( 1 \right){X_n}\left( 4 \right) + 20{X_n}\left( 2 \right){X_n}\left( 3 \right) + 24{X_n}\left( 5 \right) \\
 \end{array} \right\},\\
 &{B_6}\left( n \right) = \frac{1}{{6!}}\left\{ \begin{array}{l}
 X_n^6\left( 1 \right) + 15X_n^4\left( 1 \right){X_n}\left( 2 \right) + 40X_n^3\left( 1 \right){X_n}\left( 3 \right) + 90X_n^2\left( 1 \right){X_n}\left( 4 \right) \\
  + 144{X_n}\left( 1 \right){X_n}\left( 5 \right) + 45X_n^2\left( 1 \right)X_n^2\left( 2 \right) + 120{X_n}\left( 1 \right){X_n}\left( 2 \right){X_n}\left( 3 \right) \\
  + 40X_n^2\left( 3 \right) + 15X_n^3\left( 2 \right) + 90{X_n}\left( 2 \right){X_n}\left( 4 \right) + 120{X_n}\left( 6 \right) \\
 \end{array} \right\},\\
 &{B_m}\left( n \right) = \frac{1}{m}\sum\limits_{i = 0}^{m - 1} {{B_i}\left( n \right){X_n}\left( {m - i} \right)} .\tag{4.2}
 \end{align*}
Furthermore, we are able to obtain relations involving multiple zeta-star values and zeta values. For example, taking ${x_i} = \frac{1}{{{i^p}}},p>1\;\left( {i = 1, \cdots ,n} \right)$ in (4.2) and letting $n \to \infty $, we have
\[\zeta^\star \left( {{{\left\{ p \right\}}_m}} \right) = \frac{1}{m}\sum\limits_{i = 0}^{m - 1} {\zeta^\star \left( {{{\left\{ p \right\}}_i}} \right)\zeta \left( {pm - pi} \right)}.\tag{4.3}\]
where ${{{\left\{ p \right\}}_m}}$ denotes the $m$-tuple $\{p,...,p\}$, $\zeta^\star \left( {{{\left\{ p \right\}}_0}} \right) = 1$.
By a similar argument as in the proof of Theorem 4.1, we can get the following Theorem.
\begin{thm}
Define two sequences ${{\bar A}_m(n)}$ and ${{\bar B}_m(n)}$ by
$${\bar A}_m(n) = (m-1)!(-1)^{m-1}\underset{i = 0}{\overset{m - 1}{\sum}}(-1)^{i}\dfrac{{\bar A}_i(n)}{i!}\underset{k = 1}{\overset{n}{\sum}}x_k^{m - i}, {\bar A}_0(n)=1,$$
$${\bar B}_m(n) = \underset{k_1 = 1}{\overset{n}{\sum}}x_{k_1}\underset{k_2 = 1}{\overset{k_1-1}{\sum}}x_{k_2}\cdots\underset{k_m = 1}{\overset{k_{m-1}-1}{\sum}}x_{k_m},{\bar B}_0(n) = 1.$$
Then \[{\bar A}_m(n) = m!{\bar B}_m(n).\tag{4.4}\]
\end{thm}
\pf By induction, we can tell that ${{\bar A}_m(n)}$ and ${{\bar B}_m(n)}$ are polynomials of degree $m$ with $n$ variables $x_1, x_2, \cdots x_n$, moreover, ${\bar B}_m(n)$ has coefficient $1$ at $x_1^{s_1}x_2^{s_2}\cdots x_n^{s_n}$, when $\underset{k = 1}{\overset{n}{\sum}}s_k = m$, and $0 \leq s_k \leq 1$, otherwise, the coefficient is $0$. Suppose $c_{s_1, \cdots, s_n}$ is the coefficient of ${\bar A}_m(n)$ at $x_1^{s_1}x_2^{s_2}\cdots x_n^{s_n}$, then in order to get the conclusion, all we need to do is to prove that $c_{s_1, \cdots, s_n} = m!$ if $0 \leq s_k \leq 1$, and $c_{s_1, \cdots, s_n} = 0$ if there is a $k$ such that $s_k \geq 2$ by induction.
For convenience, we can suppose $c_{t_1, \cdots, t_n} = 0$ if there is a $t_i < 0$, otherwise $c_{t_1, \cdots, t_n}$ is the coefficient of ${\bar A}_{t_1 +\cdots + t_n}(n)$ at $x_1^{t_1}x_2^{t_2}\cdots x_n^{t_n}$.\\
For $m = 0$, the conclusion holds.
If the conclusion holds for any $i \leq m - 1$. For $m$, by the recursion formula of ${\bar A}_m(n)$, $c_{s_1, \cdots, s_n} = (m-1)!(-1)^{m-1}\underset{i = 1}{\overset{m - 1}{\sum}}\underset{k = 1}{\overset{n}{\sum}}(-1)^i\dfrac{c_{s_1,\cdots,s_k - (m - i), \cdots, s_n}}{i!}$, since for $1 \leq i \leq m-1$,  $\dfrac{{\bar A}_i(n)}{i!}\underset{k = 1}{\overset{n}{\sum}}x_k^{m - i}$ at $x_1^{s_1}x_2^{s_2}\cdots x_n^{s_n} = x_1^{s_1}\cdots x_k^{s_k - (m - i)}\cdots x_n^{s_n} \cdot x_k^{m - i}$, where $1 \leq k \leq n$, provides coefficient $\dfrac{c_{s_1,\cdots,s_k - (m - i), \cdots, s_n}}{i!}$. So,
if there are $k_1, k_2$ such that $s_{k_1}, s_{k_2} \geq 2$, then $c_{s_1,\cdots,s_k - (m - i), \cdots, s_n} = 0$ for any $0 \leq i \leq m-1, 0 \leq k \leq n$, so $c_{s_1, \cdots, s_n} = 0$.
If there is only one $k$ such that $s_k \geq 2$, say $j$, then $c_{s_1,\cdots,s_k - (m - i), \cdots, s_n} = 0$ for any $0 \leq i \leq m-1, 0 \leq k \leq n$ and $k \not = j$, so
\begin{align*}
{c_{{s_1}, \cdots ,{s_n}}} &= (m - 1)!{( - 1)^{m - 1}}\sum\limits_{i = 1}^{m - 1} {\sum\limits_{k = 1}^n {{{( - 1)}^i}\frac{{{c_{{s_1}, \cdots ,{s_k} - (m - i), \cdots ,{s_n}}}}}{{i!}}} } \\
& = (m - 1)!{( - 1)^{m - 1}}\sum\limits_{k = 1}^n {\sum\limits_{{i_k} = m - {s_k}}^{m - 1} {{{( - 1)}^{{i_k}}}\frac{{{c_{{s_1}, \cdots ,{s_k} - (m - {i_k}), \cdots ,{s_n}}}}}{{i_k!}}} } \\
& = (m - 1)!{( - 1)^{m - 1}}\sum\limits_{i = m - {s_j}}^{m - 1} {{{( - 1)}^i}\frac{{{c_{{s_1}, \cdots ,{s_j} - (m - i), \cdots ,{s_n}}}}}{{i!}}} \\
& = (m - 1)!{( - 1)^{m - 1}}({( - 1)^{m - {s_j}}}{c_{{s_1}, \cdots ,{s_{j - 1}},0,{s_{j + 1}} \cdots ,{s_n}}} + {( - 1)^{m - {s_j} + 1}}{c_{{s_1}, \cdots ,{s_{j - 1}},1,{s_{j + 1}}, \cdots ,{s_n}}})\\
& = (m - 1)!{( - 1)^{m - 1}}({( - 1)^{m - {s_j}}} + {( - 1)^{m - {s_j} + 1}})\\
& = 0.
\end{align*}
If $0\leq s_k \leq 1$ for each $1 \leq k \leq n$, then
\begin{align*}
{c_{{s_1}, \cdots ,{s_n}}} &= \left( {m - 1} \right)!{\left( { - 1} \right)^{m - 1}}\sum\limits_{i = 1}^{m - 1} {\sum\limits_{k = 1}^n {{{( - 1)}^i}\frac{{{c_{{s_1}, \cdots ,{s_k} - (m - i), \cdots ,{s_n}}}}}{{i!}}} } \\
& = \left( {m - 1} \right)!{\left( { - 1} \right)^{m - 1}}\sum\limits_{k = 1}^n {\sum\limits_{{i_k} = m - {s_k}}^{m - 1} {{{( - 1)}^{{i_k}}}\frac{{{c_{{s_1}, \cdots ,{s_k} - (m - {i_k}), \cdots ,{s_n}}}}}{{{i_k}!}}} } \\
& = \left( {m - 1} \right)!{\left( { - 1} \right)^{m - 1}}\sum\limits_{\mathop {1 \le k \le n}\limits_{{s_k} = 1} } {\sum\limits_{{i_k} = m - 1}^{m - 1} {{{( - 1)}^{{i_k}}}\frac{{{c_{{s_1}, \cdots ,1 - (m - {i_k}), \cdots ,{s_n}}}}}{{{i_k}!}}} } \\
&= \left( {m - 1} \right)!{\left( { - 1} \right)^{m - 1}}\sum\limits_{\mathop {1 \le k \le n}\limits_{{s_k} = 1} } {{{( - 1)}^{m - 1}}\frac{{{c_{{s_1}, \cdots ,{s_{k - 1}},0,{s_{k + 1}}, \cdots ,{s_n}}}}}{{(m - 1)!}}} \\
& = \left( {m - 1} \right)!\sum\limits_{\mathop {1 \le k \le n}\limits_{{s_k} = 1} } 1 \\
& = m{\rm{!}},
\end{align*}
since $c_{s_1, \cdots, s_{k-1}, 0, s_{k+1}, \cdots, s_n} = (\underset{k = 1}{\overset{n}{\sum}}s_k - 1)! = (m-1)!$ by the inductive assumption. So the conclusion holds by induction.\hfill$\square$\\
Let ${X_n}\left( m \right) := \sum\limits_{i = 1}^n {x_i^m} $ in Theorem 4.2, we have
\begin{align*}
&{{\bar B}_1}\left( n \right) = {X_n}\left( 1 \right),{{\bar B}_2}\left( n \right) = \frac{{X_n^2\left( 1 \right) - {X_n}\left( 2 \right)}}{{2!}},{{\bar B}_3}\left( n \right) = \frac{{X_n^3\left( 1 \right) - 3{X_n}\left( 1 \right){X_n}\left( 2 \right) + 2{X_n}\left( 3 \right)}}{{3!}},\\
&{{\bar B}_4}\left( n \right) = \frac{{X_n^4\left( 1 \right) - 6X_n^2\left( 1 \right){X_n}\left( 2 \right) + 8{X_n}\left( 1 \right){X_n}\left( 3 \right) + 3X_n^2\left( 2 \right) - 6{X_n}\left( 4 \right)}}{{4!}},\\
&{{\bar B}_5}\left( n \right) = \frac{1}{{5!}}\left\{ \begin{array}{l}
 X_n^5\left( 1 \right) - 10X_n^3\left( 1 \right){X_n}\left( 2 \right) + 20X_n^2\left( 1 \right){X_n}\left( 3 \right) + 15{X_n}\left( 1 \right)X_n^2\left( 2 \right) \\
  - 30{X_n}\left( 1 \right){X_n}\left( 4 \right) - 20{X_n}\left( 2 \right){X_n}\left( 3 \right) + 24{X_n}\left( 5 \right) \\
 \end{array} \right\},\\
 &{{\bar B}_6}\left( n \right) = \frac{1}{{6!}}\left\{ \begin{array}{l}
 X_n^6\left( 1 \right) - 15X_n^4\left( 1 \right){X_n}\left( 2 \right) + 40X_n^3\left( 1 \right){X_n}\left( 3 \right) - 90X_n^2\left( 1 \right){X_n}\left( 4 \right) \\
  + 144{X_n}\left( 1 \right){X_n}\left( 5 \right) + 45X_n^2\left( 1 \right)X_n^2\left( 2 \right) - 120{X_n}\left( 1 \right){X_n}\left( 2 \right){X_n}\left( 3 \right) \\
  + 40X_n^2\left( 3 \right) - 15X_n^3\left( 2 \right) + 90{X_n}\left( 2 \right){X_n}\left( 4 \right) - 120{X_n}\left( 6 \right) \\
 \end{array} \right\},\\
&{{\bar B}_m}\left( n \right) = \frac{{{{\left( { - 1} \right)}^{m - 1}}}}{m}\sum\limits_{i = 0}^{m - 1} {{{\left( { - 1} \right)}^i}{{\bar B}_i}\left( n \right){X_n}\left( {m - i} \right)}.\tag{4.5}
\end{align*}
Putting ${x_i} = \frac{1}{{{i^p}}},p>1\;\left( {i = 1, \cdots ,n} \right)$ and $n\rightarrow \infty$ in Theorem 4.2, we obtain
\[\zeta \left( {{{\left\{ p \right\}}_m}} \right) = \frac{{{{\left( { - 1} \right)}^{m - 1}}}}{m}\sum\limits_{i = 0}^{m - 1} (-1)^i{\zeta \left( {{{\left\{ p \right\}}_i}} \right)\zeta \left( {pm - pi} \right)}.\tag{4.6}\]
Taking $x_k=\frac{1}{{{{\left( {k + a} \right)}^p}}}\ (p>1)$ in (4.3) and (4.5), then letting $n\rightarrow \infty$, we also obtain the results
\[\zeta \left( {{{\left\{ p \right\}}_m};a + 1} \right) = \frac{{{{\left( { - 1} \right)}^{m - 1}}}}{m}\sum\limits_{i = 0}^{m - 1} {{{\left( { - 1} \right)}^i}\zeta \left( {{{\left\{ p \right\}}_i};a + 1} \right)\zeta \left( {pm - pi,a + 1} \right)},\tag{4.7} \]
\[{\zeta ^ \star }\left( {{{\left\{ p \right\}}_m};a + 1} \right) = \frac{1}{m}\sum\limits_{i = 0}^{m - 1} {{\zeta ^ \star }\left( {{{\left\{ p \right\}}_i};a + 1} \right)\zeta \left( {pm - pi,a + 1} \right)}.\tag{4.8} \]
By using (4.7) and (4.8), we can get the following Theorem.
\begin{thm} For integers $p>1,m\in \N$ and $a,b\in \N_0$, $a\neq -1,-2,\cdots$, we have the recurrence formulas
\begin{align*}
H\left( {a,b;m,p} \right) =& {\left( { - 1} \right)^{m - 1}}\zeta \left( {pm + 1} \right) + \frac{{{{\left( { - 1} \right)}^{m - 1}}}}{m}\sum\limits_{i = 1}^{m - 1} {{{\left( { - 1} \right)}^i}\zeta \left( {pm - pi} \right)H\left( {a,b;i,p} \right)} \\
 &+ \frac{{{{\left( { - 1} \right)}^{m - 1}}}}{m}\sum\limits_{i = 1}^{m - 1} {{{\left( { - 1} \right)}^i}\left( {m - i} \right)\zeta \left( {pm - pi + 1} \right)\zeta \left( {{{\left\{ p \right\}}_i}} \right)}, \tag{4.9}\\
 {H^ \star }\left( {a,b;m,p} \right) =& \zeta \left( {pm + 1} \right) + \frac{1}{m}\sum\limits_{i = 1}^{m - 1} {\zeta \left( {pm - pi} \right){H^ \star }\left( {a,b;i,p} \right)} \\
 & + \frac{1}{m}\sum\limits_{i = 1}^{m - 1} {\left( {m - i} \right)\zeta \left( {pm - pi+1} \right){\zeta ^ \star }\left( {{{\left\{ p \right\}}_i}} \right)}. \tag{4.10}
\end{align*}
\end{thm}
\pf In (4.7) and (4.8), taking the derivative with respect to $a$ and letting $a\rightarrow 0$ gives
\begin{align*}
H\left( {a,b;m,p} \right) =&  - \frac{1}{p}{\left. {\frac{\partial }{{\partial a}}\left( {\zeta \left( {{{\left\{ p \right\}}_m};a + 1} \right)} \right)} \right|_{a = 0}}\\
& =  - \frac{{{{\left( { - 1} \right)}^{m - 1}}}}{{pm}}\sum\limits_{i = 0}^{m - 1} {{{\left( { - 1} \right)}^i}{{\left. {\frac{\partial }{{\partial a}}\left( {\zeta \left( {{{\left\{ p \right\}}_i};a + 1} \right)\zeta \left( {pm - pi,a + 1} \right)} \right)} \right|}_{a = 0}}}, \\
{H^ \star }\left( {a,b;m,p} \right) =&  - \frac{1}{p}{\left. {\frac{\partial }{{\partial a}}\left( {{\zeta ^ \star }\left( {{{\left\{ p \right\}}_m};a + 1} \right)} \right)} \right|_{a = 0}}\\
 &=  - \frac{1}{{pm}}{\sum\limits_{i = 0}^{m - 1} {\left. {\frac{\partial }{{\partial a}}\left( {{\zeta ^ \star }\left( {{{\left\{ p \right\}}_i};a + 1} \right)\zeta \left( {pm - pi,a + 1} \right)} \right)} \right|} _{a = 0}}.
\end{align*}
By a simple calculation, we obtain the desired results.\hfill$\square$\\
Therefore, from (4.3), (4.6), (4.9) and (4.10), we give the following corollary.
\begin{cor} For integers $p\in \mathbb{N} \setminus \{1\},\ m\in \N$ and $a,b\in \N_0$, then the sums
\begin{align*}
\sum\limits_{a + b = m - 1} {\zeta \left( {{{\left\{ p \right\}}_a},p + 1,{{\left\{ p \right\}}_b}} \right)}\quad{\rm and}\quad\sum\limits_{a + b = m - 1} {{\zeta ^ \star }\left( {{{\left\{ p \right\}}_a},p + 1,{{\left\{ p \right\}}_b}} \right)}
\end{align*}
can be expressed as a rational linear combination of products of Riemann zeta values.
\end{cor}
It is easy to see that for $p>1$ and $p\in \mathbb{R}$, Corollary 4.4 ia also true.
Setting ${x_k}: = \frac{{{t^k}}}{{{k^p}}}\;\left( {p > 2} \right)$ in Theorem 4.1 and Theorem 4.2, then taking the derivative with respect to $t$ and letting $t\rightarrow 1$ and $n\rightarrow \infty$, we also obtain the following description.
\begin{cor}For real $p>2$ and integers $m\in \N,a,b\in \N_0$, then the sums
\begin{align*}
\sum\limits_{a + b = m - 1} {\zeta \left( {{{\left\{ p \right\}}_a},p- 1,{{\left\{ p \right\}}_b}} \right)}\quad{\rm and}\quad\sum\limits_{a + b = m - 1} {{\zeta ^ \star }\left( {{{\left\{ p \right\}}_a},p - 1,{{\left\{ p \right\}}_b}} \right)}
\end{align*}
can be expressed as a rational linear combination of products of Riemann zeta values.
\end{cor}
{\bf Acknowledgments.} The authors would like to thank the anonymous
referee for his/her helpful comments, which improve the presentation
of the paper.
 \end{document}